\documentclass[a4paper,reqno,12pt]{amsart}

\usepackage{relsize}

\usepackage{amsfonts,amssymb, amsmath, amsthm} 
\usepackage[fixamsmath]{mathtools} 
\usepackage{bm,bbm} 

\usepackage[spacing,kerning]{microtype}
\usepackage{ellipsis} 

\usepackage[T1]{fontenc} 
\usepackage[ansinew]{inputenc} 
\usepackage{lmodern}

\usepackage[hypcap=true]{caption}

\usepackage[left=74pt, right=70pt, top=70pt, bottom=68pt, dvips, pdftex]{geometry}

\usepackage[linktocpage,bookmarksopen,bookmarksnumbered,
      pdftitle={On small values of indefinite diagonal quadratic forms at integer points in at least five variables},
      pdfauthor={Paul Buterus, Friedrich Goetze, Thomas Hille},
      pdfsubject={Irrational Quadratic Forms},
      pdfkeywords={Oppenheim Conjecture}]{hyperref}

\usepackage[shortlabels]{enumitem} 
\setitemize{leftmargin=0.85cm}
\setenumerate{leftmargin=1cm}

\usepackage{multirow}
\newcommand\Tstrut{\rule{0pt}{2.6ex}}         
\newcommand\Bstrut{\rule[-1.8ex]{0pt}{0pt}}   

\let\oldtocsection=\tocsection
\let\oldtocsubsection=\tocsubsection

\renewcommand{\tocsection}[2]{\vspace{0.5mm} \bfseries \oldtocsection{#1}{#2}}
\renewcommand{\tocsubsection}[2]{\vspace{0.5mm}\hspace{1.5em}\oldtocsubsection{#1}{#2}}

\makeatletter
\def\subsection{\@startsection{subsection}{1} \z@{.7\linespacing\@plus\linespacing}{.5\linespacing}{\normalfont\scshape\centering}}
\renewenvironment{proof}[1][\proofname:]{\par
  \pushQED{\qed}
  \normalfont \topsep6\p@\@plus6\p@\relax
  \trivlist
  \item[\hskip\labelsep
        \itshape
    #1]\ignorespaces
}{
  \popQED\endtrivlist\@endpefalse
}
\makeatother

\DeclareFontFamily{U}{mathb}{\hyphenchar\font45}
\DeclareFontShape{U}{mathb}{m}{n}{
<-6> mathb5 <6-7> mathb6 <7-8> mathb7
<8-9> mathb8 <9-10> mathb9
<10-12> mathb10 <12-> mathb12
}{}
\DeclareSymbolFont{mathb}{U}{mathb}{m}{n}
\DeclareMathSymbol{\llcurly}{\mathrel}{mathb}{"CE}
\DeclareMathSymbol{\ggcurly}{\mathrel}{mathb}{"CF}

\newtheorem{theorem}{Theorem}[section]

\newtheorem{lemma}[theorem]{Lemma}
\newtheorem{corollary}[theorem]{Corollary}
\newtheorem{proposition}[theorem]{Proposition}
\newtheorem{assumption}[theorem]{Assumption}
\newtheorem{definition}[theorem]{Definition}

\theoremstyle{remark}
\newtheorem{remark}[theorem]{Remark}
\newtheorem*{remark-non}{Remark}

\makeatletter
\def\ps@firstpage{\ps@plain
  \def\@oddfoot{\normalfont\scriptsize \hfil\rule{0pt}{20pt}\thepage\hfil
     \global\topskip\normaltopskip}%
  \let\@evenfoot\@oddfoot
  \def\@oddhead{\@serieslogo\hss}%
  \let\@evenhead\@oddhead 
}
\makeatother

\usepackage[backend=biber,style=alphabetic,maxalphanames=99,maxnames=99]{biblatex}
\newcommand{\iu}{\mathrm{i}}
\newcommand{\eps}{\varepsilon}
\newcommand{\dif}{\mathrm{d}}
\newcommand{\abs}[1]{\lvert#1\rvert}

\newcommand{\calC}{\mathcal C}

\newcommand{\calF}{\mathcal F}
\newcommand{\calG}{\mathcal G}

\newcommand{\calJ}{\mathcal J}

\newcommand{\calO}{\mathcal O}

\newcommand{\R}{\mathbb{R}}

\newcommand{\Q}{\mathbb{Q}}
\newcommand{\N}{\mathbb{N}}

\newcommand{\Z}{\mathbb{Z}}

\newcommand{\4}{\kern1pt}

\DeclareMathOperator{\sgn}{sgn}

\newcommand*\mcap{\mathbin{\mathpalette\mcapinn\relax}}
\newcommand*\mcapinn[2]{\vcenter{\hbox{$\mathsurround=0pt\ifx\displaystyle#1\textstyle\else#1\fi\bigcap$}}}
\newcommand*\mprod{\mathbin{\mathpalette\mcaprod\relax}}  
\newcommand*\mcaprod[2]{\vcenter{\hbox{$\mathsurround=0pt\ifx\displaystyle#1\textstyle\else#1\fi\prod$}}}

\newcommand\mfrac[2]{\text{\footnotesize\raisebox{.15ex}{%
\dimen0=\fontdimen8\textfont2  
\dimen2=\fontdimen11\textfont2 
\dimen4=\fontdimen8\textfont3  
$%
\fontdimen8\textfont2=.5\dimen0
\fontdimen11\textfont2=.5\dimen2
\fontdimen8\textfont3=1.1\dimen4
\dfrac{#1}{#2}$%
\fontdimen8\textfont2=\dimen0
\fontdimen11\textfont2=\dimen2
\fontdimen8\textfont3=\dimen4
}}}

\title[On small values of indefinite diagonal quadratic forms]{On small values of indefinite diagonal quadratic forms at integer points in at least five variables}

\author[P. Buterus]{Paul Buterus}
\address{Mathematisches Institut \\ Bunsenstrasse 3-5 \\ D-37073 G\"{o}ttingen \\ Germany}
\email{buterus@mathematik.uni-goettingen.de}

\author[F. G\"{o}tze]{Friedrich G\"{o}tze}
\address{Faculty of Mathematics \\ Bielefeld University \\ P.O. Box 100131 \\ D-33501 Bielefeld \\ Germany}
\email{goetze@math.uni-bielefeld.de}

\author[T. Hille]{Thomas Hille}
\address{Mathematics Department \\ Northwestern University \\ 2033 Sheridan Road \\ Evanston \\ IL 60208 \\ USA}
\email{thomas.hille@northwestern.edu}


\numberwithin{equation}{section}

\begin{document}

\begin{abstract}
  For any $\varepsilon > 0$ we derive effective estimates for the size of a non-zero integral point $m \in \mathbb{Z}^d \setminus \{0\}$ solving the Diophantine inequality $\abs{Q[m]} < \varepsilon$, where $Q[m] = q_1 m_1^2 + \ldots + q_d m_d^2$ denotes a non-singular indefinite diagonal quadratic form in $d \geq 5$ variables. In order to prove our quantitative variant of the Oppenheim conjecture, we extend an approach developed by Birch and Davenport [BD58b] to higher dimensions combined with a theorem of Schlickewei [Sch85]. The result obtained is an optimal extension of Schlickewei's result, giving bounds on small zeros of integral quadratic forms depending on the signature $(r,s)$, to diagonal forms up to a negligible growth factor.
\end{abstract}

\maketitle


\section{Introduction}
\noindent The study of the size of the least non-trivial integral solution to homogeneous quadratic Diophantine inequalities is often referred to as the \textit{quantitative Oppenheim conjecture}; it has undergone significant developments over the past twenty years, starting with the seminal results of Bentkus and G\"otze \cite{bentkus-goetze:1997} and Eskin, Margulis and Mozes \cite{Eskin-Margulis-Mozes:1998}. Still, at present the classical result of Birch and Davenport \cite{birch-davenport:1958a} provides the sharpest known bounds within the class of diagonal forms. In the present paper we consider non-singular, indefinite, diagonal quadratic forms $Q[m] := q_1 m_1^2 + \ldots + q_d m_d^2$ of signature $(r,s)$ with $d = r+s \ge 5$ variables and generalize the result of Birch and Davenport to this class: We significantly improve the explicit bounds, established by Birch and Davenport, in terms of the signature $(r,s)$ by means of Schlickewei's work \cite{Schlickewei:1985} on the size of small zeros of integral quadratic forms (see Theorem \ref{MainTheorem}). In general, we expect the size of the least solution for real coefficients to be almost as good as for integral coefficients and, in fact, the result obtained here reflects this heuristic viewpoint.

\subsection{The Result of Birch and Davenport}
\noindent The proof used here extends a method
 developed by Birch and Davenport \cite{birch-davenport:1958a}, which in turn is a refinement of the Davenport-Heilbronn circle method \cite{davenport-heilbronn:1946}. Their approach can be used to extend bounds on small zeros of integral forms to the real case: Birch and Davenport proved in the case $d=5$ (assuming that all of the real numbers $q_1,\ldots,q_d$ are of absolute value at least one) that for any $\eps >0$ the Diophantine inequality
\begin{equation}
  \label{eq:Intro:1}
 \abs{Q[m]} = \abs{q_1 m_1^2 + \ldots + q_d m_d^2}  < \eps
\end{equation}
is non-trivially solvable in integers and furthermore gave an effective estimate on the size of the least solution: For any $\delta >0$ there is a constant $C_\delta>0$, depending on $\delta$ only, and a non-trivial integral solution $m =(m_1,\ldots,m_d) \in \Z^d \setminus \{0\} $ of \eqref{eq:Intro:1} lying in the elliptic shell defined by
\begin{equation}
  \label{eq:Intro:2}
  \abs{q_1} m_1^2 + \ldots + \abs{q_d} m_d^2 \leq C_\delta \abs{q_1 \ldots q_d}^{1+\delta} \eps^{-4-\delta}.
\end{equation}
Here the weighted norm in \eqref{eq:Intro:2} is an appropriate choice because of the scaling properties with respect to $q_1,\ldots,q_d$. More importantly, the above result implies for $d \ge 5$ and arbitrarily small $\eps >0$ that there exists a non-trivial solution of \eqref{eq:Intro:1} with integral $m_1,\ldots,m_d$ all of size $\calO(\eps^{-2-\delta})$ for any fixed $\delta >0$.

Their proof relies essentially on results on small zeros of integral forms due to these authors, see \cite{birch-davenport:1958b}. For dimension $d=5$ these bounds are in general optimal, as noted in Remark \ref{intro:remark:1}. Actually, the main ideas for proving such bounds are due to Cassels \cite{Cassels:1955}, while the modifications in \cite{birch-davenport:1958b} were of the same form as \eqref{eq:Intro:2} with the choice $\eps =1$ up to the additional dependency on $\delta$. Indeed, they showed that any indefinite quadratic form $F[m] = f_1 m_1^2 + \ldots + f_d m_d^2$ in $d \ge 5$ variables with non-zero integers $f_1, \dots, f_d$ admits a non-trivial lattice point $m=(m_1,\dots, m_d) \in \Z^d \setminus \{0\}$ satisfying
\begin{equation}
    \label{extracasselsbound}
    \begin{gathered}      
    F[m] = f_1 m_1^2 + \ldots + f_d m_d^2 =0 \quad \text{ and } \\ 0 < \abs{f_1} m_1^2 + \ldots + \abs{f_d} m_d^2 \ll_d \abs{f_1 \ldots f_d},
    \end{gathered}
\end{equation}
where we use Vinogradov's notation $\ll$ as usual. We also note here that the condition $d \ge 5$ on the dimension cannot be relaxed, since it is well-known that integral quadratic forms in four variables may fail to have non-trivial zeros.

To extend the bound \eqref{extracasselsbound} to the real case, Birch and Davenport \cite{birch-davenport:1958a}, roughly speaking, analyze regular patterns in the frequency picture of the associated counting problem by establishing rigidity in form of what we call \textit{coupled} Diophantine approximation (see Definition \ref{def:coupling} for the precise meaning) and deduce a contradiction by counting these points (i.e.\ establishing an upper and a lower bound for the number of certain Diophantine approximants) under the assumption that there are no solutions of $\abs{Q[m]} < \eps$ in the elliptic shell defined by \eqref{eq:Intro:2}.

A major feature of this approach is to avoid the evaluation of the precise sizes of the Diophantine approximants and of the absolute values of  typical quadratic Weyl sums, which are related to the approximation error via a refined Weyl inequality. In fact, an approach, which only aims for an asymptotic approximation of the number of integral solutions of \eqref{eq:Intro:1} with a sufficiently small remainder term, is not suitable for our purpose.

\subsection{Main Result: Our Extension of Schlickewei's Bound}
\noindent In view of Schlickewei's work \cite{Schlickewei:1985} on the magnitude of small zeros of integral quadratic forms - which will be the main ingredient to bound the size of the least non-trivial solution of \eqref{eq:Intro:1} - it is reasonable to expect that the exponent in the bound \eqref{eq:Intro:2} can be improved in terms of the signature $(r,s)$ and, in fact, the main objective of the present paper is to prove this extension of the work \cite{birch-davenport:1958a}.

Although the general strategy of the proof uses the approach of Birch and Davenport \cite{birch-davenport:1958a} as well, their technique fails without further analysis of certain arcs, where the Weyl sums under consideration are large, if one wishes to replace \eqref{extracasselsbound} by a better bound. To overcome this issue, we shall prove (conditionally) improved mean value estimates for certain products of Weyl sums and iterate the \textit{coupling argument} of Birch and Davenport, as we will describe in detail later (see Subsection \ref{sec:description}).

Compared to the work \cite{Cassels:1955} of Cassels, Schlickewei \cite{Schlickewei:1985} showed that the dimension, say $d_0$, of a maximal rational isotropic subspace influences the size of possible solutions essentially, rather than the mere indefiniteness (i.e.\ $d_0 \ge 1$). Moreover, by using an induction argument combined with Meyer's theorem \cite{Meyer:1884} Schlickewei derived a lower bound for $d_0$ in terms of the signature $(r,s)$ as well (see Proposition \ref{theorem:Schlickewei:2}). In Section \ref{sec:rational-case} we combine both steps of Schlickewei's work to deduce the following generalization of \eqref{extracasselsbound}. To state our modified variant of Schlickewei's result, we have to assume (w.l.o.g.) that $r \geq s$ (one can replace all $q_i$ by $-q_i$).

\begin{theorem}[Schlickewei \cite{Schlickewei:1985}]
  \label{th:Intro:SchlickeweiTheorem}
  For any non-zero integers $f_1,\ldots,f_d$, of which $r \geq 1$ are positive and $s \geq 1$ negative with $d = r+s \ge 5$, there exist integers $m_1,\ldots,m_d$, not all zero, such that
  \begin{equation}
    \label{extrabound}
    f_1 m_1^2 + \ldots + f_d m_d^2 =0 \quad \quad \text{and} \quad \quad 0 < \abs{f_1} m_1^2 + \ldots + \abs{f_d} m_d^2 \ll_d \abs{f_1 \ldots f_d}^\frac{2\beta+1}{d},
  \end{equation}
  where the exponent $\beta$ is given by
  \begin{equation}
    \label{eq:Intro:Beta}
    \beta = \beta(r,s) = \begin{cases}
      \frac{1}{2} \frac{r}{s}     & \text{for} \ r \geq s+3                \\
      \frac{1}{2} \frac{s+2}{s-1} & \text{for} \ r=s+2 \ \text{or} \ r=s+1 \\
      \frac{1}{2} \frac{s+1}{s-2} & \text{for} \ r=s
    \end{cases}
  \end{equation}
  and the implicit constant in \eqref{extrabound} depends on the dimension $d$ only.
\end{theorem}

\begin{remark}
  Compared to \eqref{extracasselsbound}, the exponent in \eqref{extrabound} is smaller for a wide range of signatures $(r,s)$ and in the cases, where the exponent is larger, we can restrict $Q$ by setting some coordinates to zero to arrive at least at the result of the case $d=5$. For example, if one has $r \sim s$, then $2 \beta \sim 1$ and therefore the exponent in \eqref{extrabound} is of order $\sim 2/d$.
\end{remark}

To simplify the analysis of \eqref{eq:Intro:1}, we may assume that $\eps=1$. Indeed, replacing all coefficients $q_j$ by $q_j/\eps$ it is sufficient to consider the solvability of the inequality
\begin{equation}
  \label{eq:Oppenheim:SmallValue}
  \abs{q_1 m_1^2 + \ldots + q_d m_d^2}  < 1.
\end{equation}
Guided by Theorem \ref{th:Intro:SchlickeweiTheorem}, we shall prove the following bound for the non-integral case, which is already comparable to \eqref{extrabound} up to the determinant of $\operatorname{diag}(f_1,\ldots,f_d)$ being replaced by the $d$-th power of the largest eigenvalue and an additional growth rate given by \eqref{eq:MeaningOfll}.

\begin{theorem}
  \label{MainTheorem}
  Let $q_1,\ldots, q_d$ be real numbers, of which $r \geq 1$ are positive and $s \geq 1$ negative, such that $\abs{q_i} \geq e^e$ and $d= r+s \geq 5$. Then there exist integers $m_1,\ldots,m_d$, not all zero, satisfying both \eqref{eq:Oppenheim:SmallValue} and
  \begin{equation}
    \label{eq:Theorem1}
    \abs{q_1} m_1^2 + \ldots + \abs{q_d} m_d^2 \llcurly_d (\max_{i=1,\ldots,d} \abs{q_i})^{1+2\beta},
  \end{equation}
  where $\beta$ is defined as in \eqref{eq:Intro:Beta}. Here the implicit constant depends on $d$ only and $A \llcurly B$ stands for
  \begin{equation}
    \label{eq:MeaningOfll}
    A \ll B^{1+\frac{20d^2}{ \log{\log{B}}}}.
  \end{equation}
\end{theorem}

The reader may note that the growth rate is considerably improved compared with \eqref{eq:Intro:2}, since we have $B^{1+\frac{20d^2}{ \log{\log{B}}}} \ll B^{1+\delta}$ for any $\delta >0$. This improvement is achieved by replacing the smoothing kernel (in the application of the circle method) by a faster decaying choice. Our result can now be summarized by the following corollary to the main theorem.

\begin{corollary}
   Let $q_1,\ldots, q_d$ be as in Theorem \ref{MainTheorem} and let $\eps >0$ be arbitrary. Then there exists a non-trivial solution $(m_1,\ldots,m_d) \in \Z^d \setminus \{0\}$ of $\abs{q_1 m_1^2 + \ldots + q_d m_d^2} < \eps$, whose size is of order $\|(m_1,\ldots,m_d)\| \llcurly_Q \eps^{-\beta}$.
\end{corollary}

Note that this result is an improved variant of the 
bound $\mathcal{O}(\eps^{-2-\delta})$ of Birch and Davenport \cite{birch-davenport:1958a} for higher dimensions in terms of the signature $(r,s)$.

\begin{remark}
  \label{intro:remark:1}
  In 1988 Schlickewei and Schmidt \cite{Schlickewei-Schmidt:1988} proved that Schlickewei's bound (in terms of the dimension $d_0$ of a maximal rational isotropic subspace) is qualitatively best possible. Of course, one can also ask if Schlickewei's bound in terms of the signature $(r,s)$ is best possible, as was already conjectured by Schlickewei himself in his first work \cite{Schlickewei:1985} on small zeros. For the class of integral quadratic forms (not necessarily diagonal) this is known for the cases $r \ge s + 3$ and $(3,2)$, see Schmidt \cite{Schmidt:1985}.
\end{remark}

\subsection{Related Results and Further Remarks}
\label{sec:rel_res}
\noindent In 1929 Oppenheim \cite{oppenheim:1929} conjectured that for any irrational quadratic form $Q$, i.e.\ $Q$ is not a real multiple of a rational form, in $d\geq 5$ variables the set $Q[\Z^d]$ contains values arbitrarily close to zero. The stronger version - conjecturing that it is sufficient to have $d \ge 3$ variables - is due to Davenport \cite{davenport-heilbronn:1946}. Actually the density of $Q[\Z^d]$ in $\R$ follows from that $Q$ either represents zero non-trivially or $Q[\Z^d]$ contains non-zero elements with arbitrarily small absolute values, provided that $d \ge 4$ and $Q$ is irrational, see \cite{oppenheim:1953i,oppenheim:1953ii,oppenheim:1953iii} and for instance  Section 5 in \cite{lewis:1973}. The validity of the conjecture was confirmed by Birch, Davenport and Ridout \cite{birch-davenport:1958c,Davenport-Ridout:1959,Ridout:1958} for $d \geq 21$ and conclusively answered in 1986 by Margulis \cite{Margulis:1989}, using methods of homogeneous dynamics. The first proof given by Margulis shows only the solvability of $\abs{Q[m]} < \eps$, whereas the modified variant, i.e.\ the solvability of $0< \abs{Q[m]} < \eps$ for irrational $Q$, was proven by Margulis subsequently as well (responding to a question by Borel). We refer to \cite{lewis:1973} and \cite{margulis:1997} for a complete historical overview until 1997.

\begin{remark}
  Baker and Schlickewei \cite{baker-schlickewei:1987} have already used Schlickewei's work \cite{Schlickewei:1985} in combination with the methods of Davenport and Ridout \cite{Davenport-Ridout:1959} to prove the Oppenheim conjecture (for not necessary diagonal forms) in some special cases ((i) $d=18$, $r=9$, (ii) $n=20$, $ 8 \le 11$, (iii) $d=20$, $7 \le r \le 13$).
\end{remark}

Nearly a decade later Eskin, Margulis and Mozes \cite{Eskin-Margulis-Mozes:1998,Eskin-Margulis-Mozes:2005} gave quantitative versions of these results, i.e.\ counting asymptotically the number of lattice points in fixed hyperbolic shells $\{ m \in \Z^d : a < Q[m] < b\}$ which are restricted to growing domains $r\Omega$ with $r \rightarrow \infty$. Such results are called \textit{quantitative Oppenheim conjecture} as well, but do not imply in a first instance explicit bounds on the size of the least non-trivial integral solution to homogeneous quadratic Diophantine inequalities: To show that the inequality $\abs{Q[m]} < \eps$ admits a non-trivial integer solution, whose size can be bounded, an effective error bound for the lattice remainder is needed. This investigation started with the work of Bentkus and G\"{o}tze \cite{bentkus-goetze:1997,bentkus-goetze:1999}, establishing effective bounds for $d \ge 9$ (however, in these works no explicit connections between certain Theta-series and Diophantine approximation of $Q$ were deduced) and later continued by G\"{o}tze and Margulis \cite{goetze-margulis:2013}. In an upcoming revised version \cite{buterus-goetze-hille-margulis:2019} we prove effective versions of the Oppenheim conjecture for $d \ge 5$ and non-diagonal forms, thereby deriving bounds on solutions of $\abs{Q[m]} < \eps$ as well, see \cite[Theorem 1.3]{buterus-goetze-hille-margulis:2019}. However, we cannot make use of the full strength of Schlickewei's bounds and therefore our Theorem \ref{MainTheorem} for diagonal forms, obtained in this paper, is sharper when compared to \cite{buterus-goetze-hille-margulis:2019}.

\begin{remark}
  We also note that weaker results, giving upper bounds in terms of the signature for general quadratic forms were established by Cook \cite{Cook:1983}, \cite{Cook:1984}, and Cook and Raghavan \cite{Cook-Raghavan:1984} using a diagonalization technique of Birch and Davenport.
\end{remark}

Recently Bourgain \cite{Bourgain:2016}, Athreya and Margulis \cite{Athreya-Margulis:2018}, and Ghosh and Kelmer \cite{Ghosh-Kelmer:2018} investigated generic variants of the quantitative Oppenheim conjecture. Bourgain \cite{Bourgain:2016} proved essentially optimal results for one-parameter families of diagonal ternary indefinite quadratic forms under the Lindel\"of hypothesis by using an analytic number theory approach. Compared to \cite{Bourgain:2016}, Ghosh and Kelmer consider in the paper \cite{Ghosh-Kelmer:2018} the space of all indefinite ternary quadratic forms, equipped with a natural probability measure, and they use an effective mean ergodic theorem for semisimple groups. In contrast, Athreya and Margulis \cite{Athreya-Margulis:2018} applied classical bounds of Rogers for $L^2$-norm of Siegel transforms to prove that for every $\delta >0$ and almost every $Q$ (with respect to the Lebesgue measure) with signature $(r,s)$, there exists a non-trival integer solution $m \in \Z^d$ of the Diophantine inequality $\abs{Q[m]} < \eps$ whose size is $\|m\| \ll_{\delta,Q} \eps^{-\frac{1}{d-2}+\delta}$, if $d \geq 3$. 

\subsection{Sketch of Proof: Extended Approach of Birch and Davenport}
\label{sec:description}
\noindent We follow the approach of Birch and Davenport \cite{birch-davenport:1958a}, which is a proof by contradiction and consists mainly of two parts: The first step is to pick out all integral solutions to the inequality $\abs{Q[m]} <1$ that are contained in a box of a certain size by integrating the product of all quadratic exponential sums $S_1,\ldots,S_d$, defined by
\begin{equation}
  \label{def:S_j}
  S_j(\alpha) := \sum_{P < \abs{q_j}^{1/2} m_j < 2d P} e(\alpha q_j m_j^2),
\end{equation}
with a suitable kernel $K$. Here we write as usual $e(x) = \exp(2\pi \iu x)$. Assuming that there are no integral solutions contained in the elliptic shell defined by
\begin{equation}
  \label{def:elliptic_shell}
  \abs{q_1} m_1^2 + \ldots + \abs{q_d} m_d^2 \le 4d^3 P^2,
\end{equation}
we deduce (in Lemma \ref{lemma:Cancellation}) that the real part of $\int_0^\infty S_1(\alpha) \ldots S_d(\alpha) K(\alpha) \, \dif \alpha$ vanishes, i.e.\ there are non-trivial cancellations in the product of the sums $S_1,\ldots,S_d$. To analyze this integral, we will divide the range of integration into four parts, namely
\begin{equation}
  \label{eq:decomp_int}
     0 < \alpha < \mfrac{1}{(8dP) q^{1/2}}; \ \ \ \mfrac{1}{(8dP) q^{1/2}} < \alpha < \mfrac{1}{(8dP) (q_0)^{1/2}}; \ \ \ \mfrac{1}{(8dP) (q_0)^{1/2}} < \alpha <  u(P)
\end{equation}
and $u(P) < \alpha$, where $q_0,q$ will be defined in \eqref{eq:DavBir:Shorcuts} and $u(P) = \log(P+\mathrm{e})^2$. First we show that on the first range the mass of the real part is highly concentrated. In fact, since $\alpha$ is `very small', van der Corput's lemma can be applied and shows that this part is at least as large as the volume of the restricted hyperbolic shell
\begin{equation}
  \label{def:volume_restriced}
  \{x \in \R^d \4 : \4 \abs{Q[x]}<1 \} \cap \{x \in \R^d \4 : \4 P < \abs{q_j}^{1/2} x_j < 2d P \ \text{ for all } \ j =1,\ldots,d \}.
\end{equation}
In comparison, the second and fourth range of the integral are negligible. Consequently the mass contained in the third range - which we will call $\calJ$ - has to be of the order of the volume of \eqref{def:volume_restriced} and hence the integral over $\calJ$ is `large' as well when integrating the absolute value of the product $S_1,\ldots,S_d$, see Lemma \ref{lemma:BirDav:Lemma7}. Moreover, it remains `large', even if we restrict ourselves to a subregion (called $\calF$) of $\calJ$, where all factors $S_1,\dots,S_d$ are uniformly `large' (see Corollary \ref{corollary:BirDav:MeasureF}).

The second step consists in finding an upper and a lower bound for the number $N_j$ of specific rational approximants $(x_i,y_i)$ of $q_i \alpha$ in this subregion of the integral. As in Birch and Davenport \cite{birch-davenport:1958b}, it is convenient to consider those parts of this subregion, where for each $i =1,\ldots,d$ both quantities $S_i(\alpha)$ and $y_i$ are all of the same magnitude independent of $\alpha$. This can be achieved by using a localization argument, i.e.\ we use a dyadic decomposition of $\calF$ into $\ll \log(P)^{2d}$ parts. In particular, we can restrict ourselves to one of these sets, say $\calG$, where the integral over $\calG$ remains `large' (see Lemma \ref{lemma:BirDav:Lemma11}).

The lower bound for $N_i$ can be derived by a standard application of the refined Weyl inequality used here, see Corollary \ref{lemma:BirDav:Lemma12}. To establish an upper bound, we shall prove on $\calG$ that $d-k$ fractions $x_i y_d/y_i x_d$ are independent of $\alpha$ (see Lemma \ref{lemma:BirDav:Lemma15}), where $k \in \{0,1,2,3\}$ depends on the size of $\beta$ and the order of magnitude of $S_{k+1},\dots, S_{d}$ (prior to that, we have already rearranged $S_1,\ldots,S_d$ in a certain way, compare \eqref{ordering_T_j}). Here $S_{k+1},\dots,S_{d}$ show a rigid behaviour as in the rational case. Indeed, the previous observation gives rise to a factorization of $x_i$ and $y_i$ as
\begin{equation}
  \label{eq:intro:factorization}
  x_i = x x_i' \quad \quad \text{and} \quad \quad y_i = y y_i'
\end{equation}
such that $x_i', y_i'$ divide a fixed number $L$, which is independent of $\alpha$. In such situations (i.e.\ if such a factorization exists) we say that $S_{k+1},\ldots,S_d$ are \textit{coupled} on $\calG$, see Definition \ref{def:coupling}.

The case $k=0$ corresponds to Birch and Davenport's paper \cite{birch-davenport:1958a}. However, this setting occurs only if $\beta \ge 2$, i.e.\ the exponent in the bound \eqref{eq:Intro:2} has to be relatively large. In fact, the main difficulty in the present paper is to overcome this issue: In Section \ref{sec:iteration} such factorizations are used to show that all pairs $(x,y)$ lie in a certain bounded set. The size of this bounded set will be substantially influenced by the exponent $\beta$; exactly at this point we are going to apply Schlickewei's bound to the integral form $x_{k+1}' y_{k+1}' m_{k+1}^2 + \ldots + x_{d}' y_{d}' m_{d}^2$, see Lemma \ref{lemma:bound_for_xy} for more details (here the factorization \eqref{eq:intro:factorization} allows us to factor out $a/q$).

As a consequence, we deduce an upper bound for the number of distinct pairs $(x,y)$, see Corollary \ref{cor:bound_fundamental}. Based on this, we establish an improved mean value estimate for $S_{k+1} \dots S_{d}$ on $\calG$, which implies better estimates for the order of magnitude of $S_{k}$. By using this improved lower bound on $S_{k}$ we will conclude that $S_{k},\ldots,S_d$ are coupled on $\calG$ as well. Now, depending on $k \in \{0,1,2,3\}$, we can iterate this argument until $k=0$ to prove that all remaining coordinates are coupled. In the course of this, we are faced with the tedious problem of comparing Schlickewei's exponent \eqref{eq:Intro:Beta} for $Q$ and all possible restrictions of $Q$ to certain subspaces with $k$ zero coordinates. This results in the number of cases listed in the \hyperref[appendix-a]{Appendix A}. To complete the proof, we deduce an inconsistent inequality (as in Birch and Davenport \cite{birch-davenport:1958a}) by establishing an upper bound for a particular $N_i$, which contradicts the lower bound found previously.

\section{Fourier Analysis and Moment Estimates}
\noindent Throughout the paper $q_1,\ldots,q_d$ denote real non-zero numbers, of which $r \ge 1$ are positive and $s \ge 1$ negative. We also introduce the notation
\begin{equation}
  \label{eq:DavBir:Shorcuts}
  q_0 = \min_{j=1,\ldots,d} \abs{q_j}, \quad \quad  q = \max_{j=1,\ldots,d} \abs{q_j} \quad \quad \text{and} \quad \quad  \abs{Q} = \mprod_{j=1}^d \abs{q_j}.
\end{equation}
Moreover, the constants throughout the proofs involved in the notation $\ll$ will not be always mentioned explicitly; these will depend on $d$ only unless stated otherwise. We also stress the underlying assumption that $d = r+s \ge 5$, since our argument depends on the solvability of non-degenerate, integral indefinite quadratic forms that are `close' to scalar multiples of $Q$. We shall ultimately deduce a contradiction from the following assumption.

\begin{assumption}
  \label{assumption:Theorem1}
  Let $q_1,\ldots, q_d$ be as introduced in Theorem \ref{MainTheorem}. Suppose that for $C_d >0$ the inequality
  \begin{equation*}
    \abs{q_1 m_1^2 + \ldots + q_d m_d^2}  < 1
  \end{equation*}
  has no solutions in integers $m_1,\ldots,m_d$, not all zero, satisfying
  \begin{equation}
    \label{eq:DomainForIntegers}
    \abs{q_1} m_1^2+ \ldots + \abs{q_d} m_d^2 \leq 4 d^3 P^2,
  \end{equation}
  where
  \begin{equation}
    \label{eq:ChoiceOfP}
    P = \exp \big\{ \big(1+\tfrac{10d^2}{\log \log H} \big) \log{H} \big\} \quad \text{and} \quad H=C_d \4 q^{\frac{1}{2}+\beta}
  \end{equation}
  and $\beta$ is defined as in \eqref{eq:Intro:Beta}.
\end{assumption}

\begin{remark}
   During the proof we will assume that the constant $C_d >0$ in Assumption \ref{assumption:Theorem1} is chosen sufficiently large. This will guarantee that the error terms under consideration are smaller (in terms of $P$, resp.\ $H$) than the leading term.
\end{remark}

  In this paper, we shall fix from now on a smoothing kernel $K := \widehat{\psi}$ with decay rate
\begin{equation}
  \label{eq:decay-rate-psi}
  \abs{\widehat{\psi}(\alpha)} \ll \exp(-\alpha/\log(\alpha+\mathrm{e})^2),
\end{equation}
where $\psi$ is a smooth symmetric probability density supported in $[-1,1]$. Compared to \cite{birch-davenport:1958a} our choice of $K$ allows to achieve the growth rate of the bound \eqref{eq:Theorem1}, since we replace the kernel by a faster decaying one. Note that the existence of such a function $\psi$ is guaranteed by the following Lemma \ref{lemma:Improvement:BirDav:Lemma1} with the choice 
\begin{equation}
 \label{eq:choice_u}
 u(\alpha) := \log(\alpha+\mathrm{e})^2.
\end{equation}

\begin{lemma}
  \label{lemma:Improvement:BirDav:Lemma1}
  Let $u$ be a positive, continuous, strictly increasing function such that
  \begin{equation}
    \label{eq:ingham:condition}
    \int_1^\infty \frac{1}{\alpha u(\alpha)} \, \dif \alpha < \infty.
  \end{equation}
  Then there exists a smooth symmetric probability density $\psi \colon \R \rightarrow [0,\infty)$ such that
  \begin{enumerate}
    \item $\psi$ supported in $[-1,1]$ and $\psi(0) \geq 1/2$,
    \item $\psi$ is increasing for $\alpha <0$ and $\psi$ decreasing for $\alpha >0$,
    \item $\abs{\widehat{\psi}(\alpha)} \ll \exp(-\abs{\alpha} u(\abs{\alpha})^{-1})$ and $\widehat{\psi}$ is real-valued and symmetric.
  \end{enumerate}
\end{lemma}

The existence of such kernels is discussed in \cite{Bhattacharya-RangaRao:2010}, see Section 10 of Chapter 2 and particularly Theorem 10.2. However, our variant cannot be found in the literature and therefore we have included a proof in \hyperref[appendix-b]{Appendix B}.

\begin{remark}
 The construction of such kernels is due to Ingham \cite{Ingham:1934} and extends the commonly used ones in the context of the circle method (compare with Lemma 1 in \cite{davenport:1956} or \cite{Bruedern-Kumchev:2001}) by using convergent infinite convolution products (instead of finitely many). As a side note, we mention that Ingham also showed that the condition \eqref{eq:ingham:condition} is necessary for the existence of such kernels.
\end{remark}

\subsection{Counting via Integration}

\noindent The starting point of Birch and Davenport's approach is the following observation.

\begin{lemma}
  \label{lemma:Cancellation}
  Assumption \ref{assumption:Theorem1} implies
  \begin{equation}
    \label{eq:Cancellation}
    \operatorname{Re} \int_0^\infty S_1(\alpha) \ldots S_d(\alpha) K(\alpha) \, \dif \alpha = 0.
  \end{equation}
\end{lemma}

\begin{proof}
  Expanding the product shows that
  \begin{equation*}
    \operatorname{Re} \int_0^\infty S_1(\alpha) \ldots S_d(\alpha) K(\alpha) \, \dif \alpha = \frac{1}{2} \sum_{P < \abs{q_1}^{\frac{1}{2}} m_1 < 2d P} \ldots \sum_{P < \abs{q_d}^{\frac{1}{2}} m_d < 2d P} \psi(q_1 m_1^2 + \ldots + q_d m_d^2).
  \end{equation*}
  Since the domain of summation is contained in \eqref{eq:DomainForIntegers}, we have
  \begin{equation*}
    \abs{q_1 m_1^2 + \ldots + q_d m_d^2} \geq 1
  \end{equation*}
  by Assumption \ref{assumption:Theorem1}. Thus, the sum is zero because $\psi$ is supported in $[-1,1]$.
\end{proof}

We begin by investigating the first range in \eqref{eq:decomp_int}, where van der Corput's lemma can be applied in order to replace the exponential sums $S_1,\ldots,S_d$ within a small part of the integration domain by analogous exponential integrals.

\begin{lemma}
  \label{lemma:BirDav:Lemma2}
  If
  \begin{equation}
    \label{eq:ConditionOnAlpha}
    0 < \alpha < (8dP)^{-1} \abs{q_j}^{-1/2},
  \end{equation}
  then we have
  \begin{equation}
    \label{eq:ApproximationIntegral}
    S_j(\alpha) = \abs{q_j}^{-1/2} I(\pm \alpha) + \calO(1),
  \end{equation}
  where the $\pm$ sign is the sign of $q_j$ and
  \begin{equation}
    I(\alpha) = \int_P^{2dP} e(\alpha \xi^2) \, \dif \xi.
  \end{equation}
\end{lemma}

\begin{proof}
  Let $f(x) = \alpha \abs{q_j} x^2$. If $P < \abs{q_j}^{1/2} x < 2d P$, then we have $f''(x) > 0$ and $0 < f'(x) < 1/2$. Hence, we can apply van der Corput's Lemma (\cite{Vinogradov:1954}, Chapter 1, Lemma 13) to get
  \begin{equation*}
    S_j(\alpha) = \int_{P\abs{q_j}^{-\frac{1}{2}}}^{2dP\abs{q_j}^{-\frac{1}{2}}} e(\alpha q_j \xi^2) \, \dif \xi +\calO(1).
  \end{equation*}
  Changing the variables of integration proves \eqref{eq:ApproximationIntegral}.
\end{proof}

\begin{lemma}
  \label{lemma:BirDav:Lemma3}
  For $\alpha > 0$ we have
  \begin{equation}
    \label{lemma:BirDav:Lemma3:eq}
    \abs{I(\pm \alpha)} \ll \min(P,P^{-1} \alpha^{-1}).
  \end{equation}
\end{lemma}

\begin{proof}
  This follows by an application of the second mean value theorem, see Lemma 3 in \cite{birch-davenport:1958a}.
\end{proof}

The next lemma, which is a generalization of Lemma 4 in \cite{birch-davenport:1958a} to dimensions greater than five, gives an upper bound for the main integral in a small neighborhood of zero.

\begin{lemma}
  \label{lemma:BirDav:Lemma4}
  We have
  \begin{equation}
    \label{lemma:BirDav:Lemma4:eq}
    \operatorname{Re} \int_0^{ (8dP)^{-1} q^{-\frac{1}{2}} } S_1(\alpha) \ldots S_d(\alpha) K(\alpha) \, \dif \alpha = M_1+R_1,
  \end{equation}
  where the main term satisfies
  \begin{equation}
    \label{eq:BirDav:MainTerm}
    M_1 \gg \delta P^{d-2} \abs{Q}^{-1/2}
  \end{equation}
  for some $\delta >0$ depending on the kernel $K$ only and the error term is bounded by
  \begin{equation}
    \label{eq:BirDav:Error1}
    \abs{R_1} \ll P^{d-3} q^{1/2} \abs{Q}^{-1/2}.
  \end{equation}
\end{lemma}

\begin{proof}
  For each $j{=}{1,\ldots,d}$ we can apply Lemma \ref{lemma:BirDav:Lemma2} in the domain of integration giving $S_j(\alpha) = \abs{q_j}^{-\frac{1}{2}} I(\pm \alpha) + \calO(1)$.
  This together with \eqref{lemma:BirDav:Lemma3:eq} of Lemma \ref{lemma:BirDav:Lemma3} yields
  \begin{equation*}
    S_j(\alpha) \ll \abs{q_j}^{-\frac{1}{2}} \min(P,P^{-1} \alpha^{-1}).
  \end{equation*}
  Thus, the error for replacing the product of all exponential sums $S_j(\alpha)$ by the product of all $\abs{q_j}^{-\frac{1}{2}} I(\pm \alpha)$ is
  \begin{equation*}
    \bigg| \prod_{j=1}^d S_j(\alpha) - \abs{q_1 \ldots q_d}^{-\frac{1}{2}} \prod_{j=1}^d I(\pm \alpha) \bigg|
    \ll \sum_{j=1}^{d-1} \sum_{\{i_1,\ldots,i_j\} \subset \{1,\ldots,d\}}   \frac{\min(P^j,P^{-j} \alpha^{-j})}{\abs{q_{i_1} \ldots q_{i_j}}^{\frac{1}{2}}}.
  \end{equation*}
  Since $\min(P,\alpha^{-1} P^{-1}) > q^{1/2}$, the right hand side is bounded by
  \begin{equation*}
    \ll q^\frac{1}{2} \abs{Q}^{-\frac{1}{2}}  \min \big( P^{d-1},P^{-(d-1)} \alpha^{-(d-1)} \big).
  \end{equation*}
  Now, up to a small error, we can replace the sum by an integral and obtain
  \begin{equation}
    \label{proof:BirDav:Lemma4:eq1}
    \begin{aligned}
      \int_0^{ (8dP)^{-1} q^{-\frac{1}{2}}} & S_1(\alpha) \ldots S_d(\alpha) K(\alpha) \, \dif \alpha \\
    &= \abs{Q}^{-\frac{1}{2}} \int_0^{ (8dP)^{-1} q^{-\frac{1}{2}}} I(\pm\alpha) \ldots I(\pm\alpha) K(\alpha) \, \dif \alpha + \calO(\Xi)
    \end{aligned}
  \end{equation}
  where
  \begin{equation*}
      \Xi :=  q^\frac{1}{2} \abs{Q}^{-\frac{1}{2}} \int_0^\infty \min(P^{d-1},P^{-(d-1)} \alpha^{-(d-1)}) \, \dif \alpha.
  \end{equation*}
  Note that the last error can be absorbed in $R_1$ by \eqref{eq:BirDav:Error1}, because it is bounded by
  \begin{equation*}
    q^\frac{1}{2} \abs{Q}^{-\frac{1}{2}} \bigg( \int_0^{P^{-2}} P^{d-1} \, \dif \alpha + \int_{P^{-2}}^\infty P^{1-d} \alpha^{1-d} \, \dif \alpha \bigg) \ll q^\frac{1}{2} \abs{Q}^{-\frac{1}{2}} P^{d-3}.
  \end{equation*}
  We can also extend the integration domain (of the integral on the right-hand side of \eqref{proof:BirDav:Lemma4:eq1}) to $\infty$, since the additional error is given by
  \begin{align*}
    \abs{Q}^{-\frac{1}{2}} \int_{ (8dP)^{-1} q^{-\frac{1}{2}}}^\infty I(\pm\alpha) \ldots I(\pm\alpha) K(\alpha) \, \dif \alpha
     & \ll \abs{Q}^{-\frac{1}{2}} \int_{ (8dP)^{-1} q^{-\frac{1}{2}}}^\infty P^{-d} \alpha^{-d} \, \dif \alpha
    \\ &\ll \abs{Q}^{-\frac{1}{2}} q^{\frac{1}{2}} P^{-1} q^{\frac{d}{2}-1}
    \ll \abs{Q}^{-\frac{1}{2}} q^{\frac{1}{2}} P^{d-3},
  \end{align*}
  where we used that $q^{1/2} < P$. Again, this error can be absorbed in $R_1$ by \eqref{eq:BirDav:Error1}.
  
  Next, we are going to establish a lower bound for the main term
  \begin{equation*}
    M_1 = \abs{q_1 \ldots q_d}^{-\frac{1}{2}} \operatorname{Re} \Big( \int_0^\infty I(\pm \alpha) \ldots I(\pm \alpha) K(\alpha) \, \dif \alpha \Big).
  \end{equation*}
  Keeping in mind that $\widehat{K}(\alpha) = \psi(-\alpha) = \psi(\alpha)$, we may rewrite the main term as
  \begin{align*}
    M_1 & = 2^{-1} \abs{Q}^{-\frac{1}{2}} \int_{P}^{2dP} \ldots \int_{P}^{2dP} \psi(\pm \xi_1^2 \pm \ldots \pm \xi_d^2) \, \dif \xi_1 \ldots \dif \xi_d                                                        \\
        & = 2^{-d-1} \abs{Q}^{-\frac{1}{2}} \int_{P^2}^{4d^2 P^2} \ldots \int_{P^2}^{4d^2 P^2} (\eta_1 \ldots \eta_d)^{-\frac{1}{2}} \psi(\pm \eta_1 \pm \ldots \pm \eta_d) \, \dif \eta_1 \ldots \dif \eta_d.
  \end{align*}
  Since $\psi(0) \geq 1/2$ (see Lemma \ref{lemma:Improvement:BirDav:Lemma1}), there exists a $\delta \in (0,1)$ such that
  \begin{equation*}
    \psi(\alpha) > 1/4 \quad \text{for all} \quad \abs{\alpha} \leq \delta.
  \end{equation*}
  Relabeling the variables, if necessary, we may suppose that the sign attached to $\eta_1$ is $+$ and that the sign attached to $\eta_2$ is $-$. As can be easily verified, the region defined by the three conditions
  \begin{equation*}
    P^2 < \eta_i < 4 P^2 \ \ \text{for} \ \ i=3,\ldots,d \quad \text{and} \quad 4(d-1) P^2 < \eta_2 < (4d(d-1)+7)P^2,
  \end{equation*}
  \begin{equation*}
    \abs{\eta_1 - \eta_2 \pm \eta_3 \pm \ldots \pm \eta_d} < \delta
  \end{equation*}
  is contained in the region of integration. Therefore, we get the lower bound
  \begin{align*}
    M_1 & > 2^{-d-3}  \abs{Q}^{-\frac{1}{2}} (2\delta) (4d^2P^2)^{-\frac{1}{2}} \int_{4(d-1) P^2}^{(4d(d-1)+7)P^2} \eta_2^{-\frac{1}{2}} \dif \eta_2 \bigg( \int_{P^2}^{4P^2} \eta^{-\frac{1}{2}} \, \dif \eta \bigg)^{d-2} \\
        & = (2^{-4} \delta) \abs{Q}^{-\frac{1}{2}} \mfrac{\sqrt{4d(d-1)+7} - \sqrt{4(d-1)}}{d} P^{d-2}
  \end{align*}
  and the latter is at least as large as $(2^{-4} \delta) \abs{Q}^{-1/2} P^{d-2}$.
\end{proof}

\subsection{Mean-Value Estimates for Quadratic Exponential Sums}

\noindent In order to guarantee that the (yet to be introduced) Diophantine approximation of $q_j \alpha$ does not vanish as well as that the resulting rational approximation of $\alpha Q$ has the same signature as $Q$, we have to extend the upper integration limit in \eqref{lemma:BirDav:Lemma4:eq} from $(8dP)^{-1} q^{-\frac{1}{2}}$ to $(8dP)^{-1} (q_0)^{-\frac{1}{2}}$. This will be done in Lemma \ref{lemma:BirDav:Lemma6} showing that the contribution of this region is (roughly) of the same order as the previous error term \eqref{eq:BirDav:Error1}.

\begin{lemma}
 \label{lemma:BirDav:Lemma6}
 We have
 \begin{equation}
  \label{eq:BirDav:Lemma6}
  R_2 = \int_{ (8dP)^{-1} q^{-\frac{1}{2}}}^{ (8dP)^{-1} (q_0)^{-\frac{1}{2}}} \abs{S_1(\alpha) \ldots S_d(\alpha)} \, \dif \alpha \ll q^\frac{1}{2} \abs{Q}^{-\frac{1}{2}} P^{d-3} (\log{P}).
 \end{equation}
\end{lemma}
  
 To prove this lemma, we will utilize both Lemma \ref{lemma:BirDav:Lemma2} and the following moment estimates for the quadratic Weyl sums $S_1,\ldots,S_d$ under consideration.
 
\begin{lemma}
 \label{lemma:BirDav:Lemma5}
 For any $n \geq 4$ we have
 \begin{equation}
   \label{eq:BirDav:Lemma5}
   \int_0^{\abs{q_j}^{-1}} \abs{S_j(\alpha)}^n \, \dif \alpha \ll \abs{q_j}^{-\frac{n}{2}} P^{n-2} (\log{P}).
 \end{equation}
\end{lemma}

\begin{proof}
 First, we use the trivial estimate $\abs{S_j(\alpha)} \ll \abs{q_j}^{-1/2}P$ to obtain
 \begin{equation*}
    \int_0^{\abs{q_j}^{-1}} \abs{S_j(\alpha)}^n \, \dif \alpha \ll \abs{q_j}^{-\frac{n-4}{2}} P^{n-4} \int_0^{\abs{q_j}^{-1}} \abs{S_j(\alpha)}^4 \, \dif \alpha
 \end{equation*}
 and subsequently we make the change of variable $\alpha = \abs{q_j}^{-1} \theta$ to get
 \begin{equation}
    \label{eq:MomentEstimate:Proof0}
    \int_0^{\abs{q_j}^{-1}} \abs{S_j(\alpha)}^n \, \dif \alpha \ll \abs{q_j}^{-\frac{n-2}{2}} P^{n-4} \int_0^1 \Big| \sum_{m \in \mathfrak{N}} e(\theta m^2) \Big|^4 \, \dif \theta,
 \end{equation}
 where the summation is taken over $\mathfrak{N} := \{m \in \N : P < \abs{q_j}^{1/2} m < 2d P \}$. Using orthogonality reveals the underlying Diophantine equation, i.e.\ the integral of the right hand side of \eqref{eq:MomentEstimate:Proof0} represents the number of solutions of
 \begin{equation}
     \label{eq:MomentEstimate:Proof1}
     v_1^2+v_2^2 = w_1^2+w_2^2,
 \end{equation}
 where $v_i,w_i \in \mathfrak{N}$ range over the interval of summation. This number can be bounded by
 \begin{equation}
      \label{eq:MomentEstimate:Proof2}
      \sum_{n <N} r^2(n)
 \end{equation}
 with $N= 8 d^2 P^2 \abs{q_j}^{-1}$. Here $r(n)$ denotes the number of representations of a natural number $n \in \N$ as a sum of two squares. As mentioned in Lemma 5 of \cite{birch-davenport:1958a}, the sum \eqref{eq:MomentEstimate:Proof2} is $\ll N \log N$. In fact, this can be proven by translating equation \eqref{eq:MomentEstimate:Proof1} into a multiplicative problem and afterwards applying the Dirichlet hyperbola method.
\end{proof}

\begin{remark} In the case $n \geq 10$ one might appeal to the Hardy-Littlewood asymptotic formula (see e.g.\ \cite{Nathanson:1996}, Theorem 5.7) and for $n \ge 6$ we could use the results in \cite{Choi-Kumchev-Osburn:2005} to drop the term $\log{N}$ as well, but this wouldn't have any effect on Theorem \ref{MainTheorem}. For completeness, we also note that the best known asymptotic formula for \eqref{eq:MomentEstimate:Proof2} can be found in \cite{Kuehleitner:1993}.
\end{remark}

A variant of our Lemma \ref{lemma:BirDav:Lemma6} is also proved in \cite{birch-davenport:1958a} under the assumption $P > \abs{Q}^{1/2}$. The situation is even easier here, since we have $P > q$. This follows directly from Assumption \ref{assumption:Theorem1} and the fact that $\beta >1/2$ or more precisely 
\begin{equation}
  \label{eq:lowerbound:beta}
  \beta \geq \frac{1}{2} \frac{d+3}{d-3} \quad \text{ if } d \text{ is odd } \quad \quad \text{and} \quad \quad \beta \geq \frac{1}{2} \frac{d+2}{d-4} \quad \text{ if } d \text{ is even.}
\end{equation}

\begin{proof}[Proof of Lemma \ref{lemma:BirDav:Lemma6}:]
  This proof does not use any properties of the quadratic form $Q$ and thus we can assume that the eigenvalues are ordered, i.e.\ $1 \leq \abs{q_1} \leq \abs{q_2} \leq \ldots \leq \abs{q_d}$. In particular, $q_0 = \abs{q_1}$ and $q = \abs{q_d}$. We begin by splitting the interval of integration into the $d-1$ intervals
  \begin{equation*}
    I_k = \{\alpha \in (0,\infty) : (8dP\abs{q_k}^\frac{1}{2})^{-1} < \alpha < (8dP \abs{q_{k-1}}^\frac{1}{2})^{-1} \},
  \end{equation*}
  where $k=2,\ldots,d$. If $j \leq k-1$, then the condition \eqref{eq:ConditionOnAlpha} of Lemma \ref{lemma:BirDav:Lemma2} is satisfied. Therefore, combined with Lemma \ref{lemma:BirDav:Lemma3}, we obtain for $\alpha \in I_k$ the inequality
  \begin{equation}
    \label{eq:BirDav:27}
    \abs{S_j(\alpha)} \ll \abs{q_j}^{-\frac{1}{2}} P^{-1} \alpha^{-1}+1 \ll \abs{q_j}^{-\frac{1}{2}} P^{-1} \alpha^{-1}.
  \end{equation}
  For $j \geq k$ we use the trivial estimate $\abs{S_j(\alpha)} \ll P \abs{q_j}^{-\frac{1}{2}}$ to conclude that
  \begin{equation*}
    \abs{S_1(\alpha) \ldots S_d(\alpha)} \ll \abs{Q}^{-\frac{1}{2}} (P\alpha)^{1-k} P^{d-(k-1)}.
  \end{equation*}
  If $k \geq 3$, then we find the bound
  \begin{equation*}
    \int_{I_k} \abs{S_1(\alpha) \ldots S_d(\alpha)} \, \dif \alpha \ll \abs{Q}^{-\frac{1}{2}} P^{d-2(k-1)} (P \abs{q_k}^\frac{1}{2})^{k-2} = \abs{Q}^{-\frac{1}{2}} P^{d-2} (P^{-1} \abs{q_k}^\frac{1}{2})^{k-2}
  \end{equation*}
  and this is smaller than $\ll \abs{Q}^{-1/2} q^{1/2} P^{d-3}$. Next we treat the case $k=2$ corresponding to the interval $I_2$. For $j=1$ the inequality \eqref{eq:BirDav:27} still holds and therefore
  \begin{equation}
    \label{eq:EstmateSum1}
    \abs{S_1(\alpha)} \ll \abs{q_1}^{-\frac{1}{2}} P^{-1} \alpha^{-1} \ll \abs{q_1}^{-\frac{1}{2}} \abs{q_2}^{\frac{1}{2}}.
  \end{equation}
  Let $j=2,\ldots,d$. Dividing $I_2$ into parts of length $\abs{q_j}^{-1}$, i.e.\ the period of $S_j$, gives
  \begin{equation*}
    \int_{I_2} \abs{S_j(\alpha)}^{d-1} \, \dif \alpha \leq (1+ \abs{q_j} (8dP \abs{q_1}^\frac{1}{2})^{-1}) \int_0^{\abs{q_j}^{-1}} \abs{S_j(\alpha)}^{d-1} \, \dif \alpha \ll \int_0^{\abs{q_j}^{-1}} \abs{S_j(\alpha)}^{d-1} \, \dif \alpha,
  \end{equation*}
  where $P \ge q$ was used. Next we apply the mean value estimates, mentioned in Lemma \ref{lemma:BirDav:Lemma5}, to deduce that 
 \begin{equation*}
  \int_{I_2} \abs{S_j(\alpha)}^{d-1} \, \dif \alpha  \ll \abs{q_j}^{-\frac{d-1}{2}} P^{d-3} (\log{P})
 \end{equation*}
 and use H\"{o}lder's inequality to obtain
 \begin{equation*}
    \int_{I_2} \abs{S_2(\alpha) \ldots S_d(\alpha)} \, \dif \alpha \ll \abs{q_2 \ldots q_d}^{-\frac{1}{2}} P^{d-3} (\log{P}).
 \end{equation*}
 Together with equation \eqref{eq:EstmateSum1} we find
 \begin{equation*}
   \int_{I_2} \abs{S_1(\alpha) \ldots S_d(\alpha)} \, \dif \alpha \ll \abs{q_2}^\frac{1}{2} \abs{Q}^{-\frac{1}{2}} P^{d-3} (\log{P}) \ll q^\frac{1}{2} \abs{Q}^{-\frac{1}{2}} P^{d-3} (\log{P}). \qedhere
 \end{equation*}
\end{proof}

We end this subsection by combining the previous estimates in order to prove

\begin{lemma}
 \label{lemma:BirDav:Lemma7}
 Under Assumption \ref{assumption:Theorem1}, we may choose $C_d \gg 1$, occurring in the definition of $P$ in \eqref{eq:ChoiceOfP}, such that
 \begin{equation}
   \label{eq:MainContributionIntegral}
   \int_{(8dP)^{-1} (q_0)^{-\frac{1}{2}} }^{u(P)} \abs{S_1(\alpha) \ldots  S_d(\alpha) K(\alpha)} \, \dif \alpha \gg \abs{Q}^{-\frac{1}{2}} P^{d-2}.
 \end{equation}
\end{lemma}

In particular, we may neglect the tail of the integral  using the decay of $K$, see \eqref{eq:decay-rate-psi} and \eqref{eq:choice_u} for the definition  of $u(P)$.

\begin{proof}
  According to Lemmas \ref{lemma:Cancellation}, \ref{lemma:BirDav:Lemma4} and \ref{lemma:BirDav:Lemma6} we have
  \begin{equation*}
    M_1+M_2+R_1+R_2+R_3=0,
  \end{equation*}
  where $M_1 \gg \abs{Q}^{-\frac{1}{2}} P^{d-2}$, $\abs{R_1}+\abs{R_2} \ll q^{\frac{1}{2}} \abs{Q}^{-\frac{1}{2}} P^{d-3} (\log{P}) \ll \abs{Q}^{-\frac{1}{2}} P^{d-\frac{5}{2}}$
  and
  \begin{align*}
    R_3 &= \operatorname{Re} \int_{u(P)}^\infty S_1(\alpha) \ldots S_d(\alpha) K(\alpha) \, \dif \alpha, \\ 
    M_2 &= \operatorname{Re} \int_{ (8dP)^{-1} (q_0)^{-\frac{1}{2}} }^{u(P)} S_1(\alpha) \ldots  S_d(\alpha) K(\alpha) \, \dif \alpha.
  \end{align*}
  We can easily bound the tail $R_3$: Using the trivial estimate $\abs{S_j(\alpha)} \ll P \abs{q_j}^{-\frac{1}{2}}$ and the decay of $K$ gives
  \begin{equation*}
    R_3 \ll P^d \abs{Q}^{-\frac{1}{2}} \int_{u(P)}^\infty \exp(-\alpha u(\alpha)^{-1}) \, \dif \alpha \ll \abs{Q}^{-\frac{1}{2}} P^{d-3}.
  \end{equation*}
  Combining the previous estimates we end up with
  \begin{equation*}
    \abs{M_1+M_2} \leq \abs{R_1}+\abs{R_2}+\abs{R_3} \ll \abs{Q}^{-\frac{1}{2}} P^{d-3} \big( 1+ P^\frac{1}{2} \big).
  \end{equation*}
  Thus, we may increase $C_d \gg 1$ such that the latter term is smaller than the lower bound for $M_1$, and conclude that
  \begin{equation*}
    P^{d-2} \abs{Q}^{-\frac{1}{2}} \ll \abs{M_2} \leq \int_{ (8dP)^{-1} (q_0)^{-\frac{1}{2}} }^{u(P)} \abs{S_1(\alpha) \ldots  S_d(\alpha) K(\alpha)} \, \dif \alpha. \qedhere
  \end{equation*}
\end{proof}

\subsection{Ordering and Contribution of the Peaks}
\noindent In the following we will show that the main contribution to the integral \eqref{eq:MainContributionIntegral} arises from a certain subregion on which every $S_1,\ldots,S_d$ is large. Before doing this, we shall fix an ordering of $S_1,\ldots, S_d$ as well, which will be necessary in order to perform the coupling argument and its iteration. For this, we define
\begin{equation}
  \label{eq:Def:calJ}
  \calJ := \{\alpha \in (0,\infty) : (8dP q_0^\frac{1}{2})^{-1} < \alpha < u(P)\}
\end{equation}
and write
\begin{equation}
  \label{eq:Def:Permutation}
  \calJ_\pi := \{\alpha \in \calJ : \abs{q_{\pi(1)}}^\frac{1}{2} \abs{S_{\pi(1)}(\alpha)} \leq \ldots \leq \abs{q_{\pi(d)}}^\frac{1}{2} \abs{S_{\pi(d)}(\alpha)} \}
 \end{equation}
for any permutation $\pi$ of the set $\{1,\ldots,d\}$. Obviously, all these sets cover $\calJ$ completely and since there are only finitely many permutations of $\{1,\ldots,d\}$ Lemma \ref{lemma:BirDav:Lemma7} implies

\begin{lemma}
  \label{lemma:existence_ordering}
  Under Assumption \ref{assumption:Theorem1}, there exists a permutation $\pi$ of the set $\{1,\ldots,d\}$ such that
  \begin{equation}
    \label{eq:Measure_Ordering}
    \int_{\calJ_\pi} \abs{S_1(\alpha) \ldots S_d(\alpha) K(\alpha)} \, \dif \alpha \gg P^{d-2} \abs{Q}^{-\frac{1}{2}}.
  \end{equation}
\end{lemma}

From now on we shall fix a permutation $\pi$ satisfying \eqref{eq:Measure_Ordering}. With this ordering at hand, we are in position to prove that the integral in \eqref{eq:Measure_Ordering} can be restricted to
\begin{equation}
  \label{eq:Def:calF}
  \calF := \{\alpha \in \calJ_\pi : \abs{q_{\pi(i)}}^\frac{1}{2} \abs{S_{\pi(i)}(\alpha)} > P (u(P)^2 q)^{-\kappa(i)} \text{ for all }  i=1,\ldots,d\},
\end{equation}
where $\kappa(i) := \min\{i,(d-4)\}^{-1}$. Indeed, we have
 
\begin{lemma}
 \label{lemma:BirDav:Lemma10}
 Independently of Assumption \ref{assumption:Theorem1}, the estimate
 \begin{equation}
   \label{eq:BirDav:Lemma10:Bound}
   \int_{\calJ_\pi \setminus \calF} \abs{S_1(\alpha) \ldots S_d(\alpha)} \, \dif \alpha \ll \abs{Q}^{-\frac{1}{2}} P^{d-2} (\log{P})^{-1}
 \end{equation}
 holds, where the error term depends on the dimension $d$ only.
\end{lemma}

\begin{proof}
  First we cover the complement $\calJ_\pi \setminus \calF$ by the (not necessarily disjoint) union of $d$ many sets given by
  \begin{equation*}
    \calC_j :=  \{\alpha \in \calJ_\pi : \abs{q_{\pi(j)}}^{\frac{1}{2}} \abs{S_{\pi(j)}(\alpha)} \leq P (u(P) q)^{- \kappa(j)} \},
  \end{equation*}
  where $j=1,\ldots,d$. If $\alpha \in \calC_j$, then \eqref{eq:Def:Permutation} implies that
  \begin{equation*}
     \abs{q_{\pi(1)}}^{1/2} \abs{S_{\pi(1)}(\alpha)} \leq \ldots \leq \abs{ q_{\pi(j)}}^{1/2} \abs{S_{\pi(j)}(\alpha)} 
  \end{equation*}
  and therefore the left hand side of \eqref{eq:BirDav:Lemma10:Bound}, restricted to the region $\calC_j$, is bounded by
  \begin{equation}
    \label{eq:proof:Lemma10:eq1}
    \ll \abs{q_{\pi(1)} \ldots q_{\pi(k)}}^{-\frac{1}{2}} P^{k} (u(P)^2 q)^{-1} \int_0^{u(P)} \abs{S_{\pi(k+1)}(\alpha) \ldots S_{\pi(d)}(\alpha)} \, \dif \alpha,
  \end{equation}
  where $k= \min(j,d-4)$. This choice of $k$ permits to apply Lemma \ref{lemma:BirDav:Lemma5}: Since $S_i$ is a periodic function with period $\abs{q_i}^{-1}$, we find
  \begin{equation*}
    \int_0^{u(P)} \abs{S_i(\alpha)}^{d-k} \, \dif \alpha \ll u(P) \abs{q_i} \int_0^{\abs{q_i}^{-1}} \abs{S_i(\alpha)}^{d-k} \, \dif \alpha \ll q P^{d-k-2} \abs{q_i}^{-(d-k)/2} u(P) (\log{P}).
  \end{equation*}
  Thus, we can make use of H\"{o}lder's inequality to obtain
  \begin{equation*}
    \int_0^{u(P)} \abs{S_{\pi(k+1)}(\alpha) \ldots S_{\pi(d)}(\alpha)} \, \dif \alpha \ll q \4 \abs{q_{\pi(k+1)} \ldots q_{\pi(d)}}^{-\frac{1}{2}} P^{d-k-2} \4 u(P) \4 (\log{P})
  \end{equation*}
  and combined with \eqref{eq:proof:Lemma10:eq1} we conclude that
  \begin{equation*}
    \int_{\calJ_\pi \setminus \calF} \abs{S_1(\alpha) \ldots S_d(\alpha)} \, \dif \alpha \ll \abs{Q}^{-\frac{1}{2}} P^{d-2} (\log{P})^{-1}. \qedhere
  \end{equation*}
\end{proof}
 
 Everything considered, applying both Lemma \ref{lemma:existence_ordering} combined with Lemma \ref{lemma:BirDav:Lemma10}, proves the following corollary.
 
 \begin{corollary}
  \label{corollary:BirDav:MeasureF}
  Under Assumption \ref{assumption:Theorem1}, we may increase the constant $C_d \gg 1$, occurring in the definition of $P$ in \eqref{eq:ChoiceOfP}, such that
  \begin{equation}
   \int_\calF \abs{S_1(\alpha) \ldots S_d(\alpha) K(\alpha) } \, \dif \alpha \gg P^{d-2} \abs{Q}^{-\frac{1}{2}}.
  \end{equation}
 \end{corollary}

\begin{remark}
 \label{remark:cannot_improve}
 We note that the usual proof of the Hardy-Littlewood asymptotic formula shows that the mean value estimates, used here for products of $S_1, \ldots,S_d$, are in general (up to $\log$ factors) best possible. In particular, one cannot improve the exponent $\kappa(i)$ without using additional information regarding the underlying quadratic form $Q[m] = q_1 m_1^2+ \ldots + q_d m_d^2$. To obtain better moment estimates (as in Lemma \ref{averagewithresonance}) we need to iterate the coupling argument of Birch and Davenport and exploit Assumption \ref{assumption:Theorem1} to 'couple' almost all coordinates (in the sense of Definition \ref{def:coupling}) and afterwards count certain arcs (see Lemma \ref{lemma:bound_for_xy}, resp.\ Corollary \ref{cor:bound_fundamental}).
\end{remark}

\section{First Coupling via Diophantine Approximation}
\label{sec:First:DioApprCoup}
\noindent As shown in Corollary \ref{corollary:BirDav:MeasureF}, the integral over $\calF$ is relatively large. Now we shall split $\calF$ further into parts, where the quantity $S_j$ has a specified order of magnitude in terms of the following Diophantine approximation: By Dirichlet's approximation theorem, applied to any $\alpha \in \calF$ and $j \in \{1,\ldots,d\}$, there exist coprime integer pairs $(x_j,y_j)$ such that
\begin{equation}
  \label{eq:Dio:1}
  q_j \alpha = \frac{x_j}{y_j} + \rho_j \quad \quad \text{and} \quad \quad 0 < y_j \leq 8d P \lvert q_j \rvert^{-\frac{1}{2}},
\end{equation}
where the approximation error is bounded by
\begin{equation}
  \label{eq:Dio:2}
  \abs{\rho_j} < y_j^{-1} (8d P \abs{q_j}^{-\frac{1}{2}})^{-1}.
\end{equation}
For convenience, we introduce the following notations as well: We shall denote by $\Z^2 _{\text{prim}}$ the set of coprime integral pairs $(x,y)$ with $y >0$ and for any $\alpha \in \R$ we define
\begin{equation*}
  \mathfrak{D}_j(\alpha) := \{ (x_j,y_j) \in \Z^2_{\text{prim}} \, : \, (x_j,y_j) \text{ are chosen as in \eqref{eq:Dio:1} satisfying \eqref{eq:Dio:2}} \}.
\end{equation*}
Note here that none of $x_1,\ldots,x_d$ are zero, since $\abs{q_j} \alpha > \abs{q_j} (8dP)^{-1} (q_0)^{-\frac{1}{2}} > \abs{\rho_j}$ holds in the integration region $\calF$ of interest and thus $\abs{x_j} \ge y_j (\abs{\alpha q_j} - \abs{\rho_j}) > 0$.
 
\subsection{Refined Variant of Weyl's Inequality}
\noindent In order to control the magnitude of $S_j(\alpha)$ in terms of the denominator $y_j$ and the approximation error $\rho_j$, corresponding to the Diophantine approximation of $q_j \alpha$ introduced previously, we need the following (well-known) variant of Weyl's inequality.

\begin{lemma}
 \label{lemma:BirDav:Corollary2}
 If \eqref{eq:Dio:1} and \eqref{eq:Dio:2} hold, then we have
 \begin{equation}
     \label{eq:lemma:BirDav:Corollary2}
     \abs{S_j(\alpha)} \ll (y_j)^{-\frac{1}{2}} (\log{P}) \min( P \abs{q_j}^{-\frac{1}{2}},P^{-1} \abs{q_j}^\frac{1}{2}\abs{\rho_j}^{-1}).
 \end{equation}
\end{lemma}
 
 Lemma \ref{lemma:BirDav:Corollary2} can be proved along  the same lines as the corollary following Lemma 9 in \cite{birch-davenport:1958a}. Nevertheless, for completeness, we have included the proof here.
 
\begin{proof}[Sketch of Proof:]
 This is a corollary of the subsequent Lemma \ref{lemma:BirDav:Lemma9}: By taking $A= P \abs{q_j}^{-1/2}$, $x=x_j$, $y =y_j$ and replacing $\alpha$ by $q_j \alpha$ we find that
 \begin{equation*}
     S_j(\alpha) = y_j^{-1} \bigg( \sum_{m=1}^{y_j} e(x_j m^2/ y_j) \bigg) \int_{P \abs{q_j}^{-1/2}}^{2d P \abs{q_j}^{-1/2}} \exp(2\pi \iu \rho_j \xi^2) \, \dif \xi + \calO(y_j^{1/2} \log 2y_j).
  \end{equation*}
  In view of \eqref{eq:Dio:1} and \eqref{eq:Dio:2} we have $y_j^{1/2} \ll y_j^{-1/2} P \abs{q_j}^{-1/2}$ and $y_j^{1/2} \ll y_j^{-1/2} P^{-1} \abs{q_j}^{1/2} \abs{\rho_j}^{-1}$. Thus, combined with $\log(2y_j) \ll \log(P)$, we infer that the $\calO$-term is negligible (compared to the right side of \eqref{eq:lemma:BirDav:Corollary2}). To complete the proof of \eqref{eq:lemma:BirDav:Corollary2}, one has to use $S_{x_j,y_j} \ll y_j^{1/2}$ (this result is well known, see e.g.\ Lemma 8 in \cite{birch-davenport:1958a}) and afterwards estimate the integral as in the proof of Lemma \ref{lemma:BirDav:Lemma3}.
\end{proof}

\begin{lemma}
 \label{lemma:BirDav:Lemma9}
 Suppose that $A \gg 1$ and that $\alpha \in \R$ is a real number satisfying
 \begin{equation}
   \alpha = \frac{x}{y} + \rho,
 \end{equation}
 where $(x,y) \in \Z^2_{\text{prim}}$ are coprime integers with
 \begin{equation}
   \label{lemma:BirDav:Lemma9:eq2}
   0 < y \ll A \quad \text{and} \quad (8d) \4 \abs{\rho} < y^{-1} A^{-1}.
 \end{equation}
 Then
 \begin{equation}
   \label{lemma:BirDav:Lemma9:eq}
   \sum_{A < m < 2dA} e(\alpha m^2) = y^{-1} \bigg( \sum_{m=1}^y e(am^2/y) \bigg) \int_{A}^{2dA} e(\rho \xi^2) \, \dif \xi + \calO(y^\frac{1}{2} \log{2y}).
 \end{equation}
\end{lemma}

We omit the proof here since it is given in Lemma 9 of \cite{birch-davenport:1958a} with the following minor changes:  The endpoints of summation and integration must be adjusted while noting that the condition $1/(2 y \abs{\rho}) > 4dA$ has to be fulfilled. But this is certainly the case because \eqref{lemma:BirDav:Lemma9:eq2}.

\begin{remark}
  The main procedure in \cite{birch-davenport:1958a}  is to split the sum on the left hand side of \eqref{lemma:BirDav:Lemma9:eq} according to the residue classes $\operatorname{mod} \, q$ and then apply Poisson's summation formula to each of these sums. A well known alternative is to use a truncated form of the Poisson summation formula, see Lemma 4.2 and Theorem 4.1 in \cite{Vaughan:1997}.
\end{remark}

\subsection{Localization of the Set \texorpdfstring{$\calF$}{F}}

\noindent Here we aim to further localize the region $\calF$ by using a dyadic decomposition  according to the size of $\abs{S_1(\alpha)},\ldots,\abs{S_d(\alpha)}$ and $y_1,\ldots,y_d$ as follows: For each $j=1,\ldots,d$ let $T_j = 2^{t(j)}$ and $U_j = 2^{u(j)}$ denote dyadic numbers with integer exponents $t(j),u(j) \in \Z$. Corresponding to these numbers we introduce the sets
\begin{equation}
 \label{eq:DefSet:G}
 \begin{aligned}
  \calG(T_1,\ldots,T_d,U_1,\ldots,U_d) = \big\{ \alpha \in \calF :  \forall j \! \in \! \{1,\ldots,d\} \, \exists (x_j,y_j) \in \mathfrak{D}_j(\alpha) \text{ such that } &\\
  T_j P/2 < \abs{q_j}^\frac{1}{2} \abs{S_j(\alpha)} \leq T_j P \text{ and } U_j/2 < y_j \leq U_j & \big\}.
 \end{aligned}
\end{equation}
In what follows we shall assume, for notational simplicity, that the coordinates are relabeled such that \eqref{eq:Def:Permutation} holds for the trivial permutation and, as a consequence, we can write
\begin{equation}
  \label{ordering_T_j}
  T_1 \ll \ldots \ll T_d.
\end{equation}
Additionally, we have only to consider those sets $\calG(T_1,\ldots,T_d,U_1,\ldots,U_d)$ which are not-empty and for any $\alpha \in  \calG(T_1,\ldots,T_d,U_1,\ldots,U_d)$ one can see that
\begin{equation}
  \label{eq:BoundsForT_j}
  (u(P)^2 q)^{-\kappa(j)} < T_j < 4d,
\end{equation}
where we used, on the one hand, the trivial upper bound $\abs{S_j(\alpha)} \leq 2dP\abs{q_j}^{-1/2}$ and, on the other hand, the lower bound in \eqref{eq:Def:calF}. Of course, we have 
\begin{equation*}
  U_j \geq y_j \geq 1,
\end{equation*}
i.e.\ $u(j) \in \N_0$. Moreover, we may apply Lemma \ref{lemma:BirDav:Corollary2} to obtain
\begin{equation*}
  T_j \ll (y_j)^{-1/2} (\log{P}) \min \big(1, P^{-2} \abs{q_j} \abs{\rho_j}^{-1} \big) \ll U_j^{-1/2} (\log{P}) \min \big(1, P^{-2} \abs{q_j} \abs{\rho_j}^{-1} \big).
\end{equation*}
Hence we find that
\begin{equation}
  \label{eq:UpperBoundForU_j}
  U_j \ll (\log{P})^2 T_j^{-2}
\end{equation}
and
\begin{equation}
   \label{eq:BirDav:eq56}
   \abs{q_j}^{-1} \abs{\rho_j} \ll P^{-2} (\log{P}) T_j^{-1} U_j^{-1/2}.
\end{equation}

\begin{lemma}
  \label{lemma:BirDav:Lemma11}
  Under Assumption \ref{assumption:Theorem1}, there exist numbers $T_1,\ldots,T_d, U_1,\ldots,U_d$ such that
  \begin{equation}
    \label{eq:BirDav:Lemma11}
    \int_{\calG(T_1,\ldots,T_d,U_1,\ldots,U_d)} \abs{S_1(\alpha) \ldots S_d (\alpha) K(\alpha)} \, \dif \alpha \gg \abs{Q}^{-\frac{1}{2}} P^{d-2} (\log{P})^{-2d}.
  \end{equation}
\end{lemma}

\begin{proof}
  On the one hand, we know from Corollary \ref{corollary:BirDav:MeasureF} that
  \begin{equation*}
    \int_{\calF} \abs{S_1(\alpha) \ldots S_d (\alpha) K(\alpha)} \, \dif \alpha \gg \abs{Q}^{-\frac{1}{2}} P^{d-2}.
  \end{equation*}
  On the other hand, \eqref{eq:BoundsForT_j} implies
  \begin{equation*}
    1 \gg t(j) \gg - \log{\log{P}} - \log{q} \gg - \log{P},
  \end{equation*}
  and combined with \eqref{eq:UpperBoundForU_j} we find
  \begin{equation*}
    0 \leq u(j) \ll \log \log P +\abs{t(j)} \ll \log{P}.
  \end{equation*}
  Hence, the minimal number of choices for $T_1,\ldots,T_d,U_1,\ldots,U_d$ to cover all $\calF$ is $\ll (\log{P})^{2d}$. In particular, there is at least one choice of $T_1,\ldots,T_d,U_1,\ldots,U_d$ with
  \begin{equation*}
    \int_{\calG(T_1,\ldots,T_d,U_1,\ldots,U_d)} \abs{S_1(\alpha) \ldots S_d (\alpha) K(\alpha)} \, \dif \alpha \gg \abs{Q}^{-\frac{1}{2}} P^{d-2}  (\log{P})^{-2d}. \qedhere
  \end{equation*}
\end{proof}

Here and subsequently, we fix a choice of such $T_1,\ldots,T_d,U_1,\ldots,U_d$, satisfying \eqref{eq:BirDav:Lemma11} of Lemma \ref{lemma:BirDav:Lemma11}, and write
\begin{equation}
  \calG=\calG(T_1,\ldots,T_d,U_1,\ldots,U_d).
\end{equation}
Moreover, for each $j \in \{1,\ldots,d\}$ let
\begin{equation*}
  N_j := \# \{ (x_j,y_j) \in \Z^2 _{\text{prim}} \, : \, \exists \alpha \in \calG \text{ such that } (x_j,y_j) \in \mathfrak{D}_j(\alpha) \}
\end{equation*}
denote the number of distinct integer pairs $(x_j,y_j) \in \mathcal{D}_j(\alpha)$ which arise from all $\alpha \in \calG$. The previous Lemma \ref{lemma:BirDav:Lemma11} leads to the next lower bound on $N_j$.

\begin{corollary}
  \label{lemma:BirDav:Lemma12}
  For any fixed numbers $T_1,\ldots,T_d$, $U_1,\ldots,U_d$, satisfying \eqref{eq:BirDav:Lemma11}, we have 
  \begin{equation}
    \label{eq:BirDav:Lemma12:Count:N_j}
    N_j \gg (\log{P})^{-2d} (T_1 \ldots T_d)^{-1} (T_j U_j^{1/2}).
  \end{equation}
\end{corollary}

\begin{proof}
  If $\alpha \in \calG(T_1,\ldots,T_d,U_1,\ldots,U_d)$, then \eqref{eq:DefSet:G} shows that
  \begin{equation*}
    \abs{S_1(\alpha) \ldots S_d(\alpha)} \ll \abs{Q}^{-1/2} P^d (T_1 \ldots T_d)
  \end{equation*}
  and therefore the bound \eqref{eq:BirDav:Lemma11} implies
  \begin{equation}
    \label{old:BirDav:Lemma11}
    \abs{\calG} \gg P^{-2} (T_1 \ldots T_d)^{-1} (\log{P})^{-2d}.
  \end{equation}
  At the same time, inequality \eqref{eq:BirDav:eq56} implies that for each integer pair $(x_i,y_j) \in \mathfrak{D}_j(\alpha)$, arising from $\alpha \in \calG$, $\alpha$ is located in an interval of length bounded by $\ll P^{-2} T_j^{-1} U_j^{-1/2}$. Thus $\abs{\calG} \ll N_j (P^{-2} T_j^{-1} U_j^{-1/2})$ and together with \eqref{old:BirDav:Lemma11} and a simple rearrangement the claimed inequality \eqref{eq:BirDav:Lemma12:Count:N_j} follows.
\end{proof}

\subsection{Coupling of the Rational Approximants}
\label{subec:coupling}
\noindent In the following we shall establish that at least $d-3$ coordinates are coupled and later on iterate this argument to deduce that all coordinates are coupled. To be precise, we define \textit{coupling} as follows.

\begin{definition}
  \label{def:coupling}
  Let $1 \! \leq \!  j_1 \! < \!  \ldots \! < \!  j_k \! \leq \! d$, where $k \in \{1,\ldots,d\}$. We say that the coordinates $j_1,\ldots,j_k$ associated to $q_{j_1},\ldots,q_{j_k}$ (resp.\ the exponential sums $S_{j_1},\ldots,S_{j_k}$) can be coupled if for any $\alpha \in \calG$ and $j \in \{j_1,\ldots, j_k\}$ the pairs $(x_j,y_j) \in \mathfrak{D}_j(\alpha) $ are of the form
  \begin{equation}
    x_j = x x_j' \quad \quad \text{and} \quad \quad y_j = y y_j',
  \end{equation}
  where $x,y>0$ are coprime integers and $x_j',y_j'$ divide some integer $L \in \N$ such that $L$ is independent of $\alpha \in \calG$.
\end{definition}

The following lemma on the number of rational approximants with bounded denominator will be the key tool for the first coupling argument and later on for its iteration as well.

\begin{lemma}
  \label{lemma:dioph_approx}
  Let $\eta >0$, $X >0$ and $\theta$ be real numbers, such that there exist $N$ distinct integer pairs $(x,y)$ satisfying
  \begin{equation}
    \abs{\theta x - y} < \eta \quad \quad \text{and} \quad \quad 0 < \abs{x} < X.
  \end{equation}
  Then either $N<24 \eta X$ or all integer pairs $(x,y)$ have the same ratio $y/x$.
\end{lemma}

\begin{proof}
  This is Lemma 14 in \cite{birch-davenport:1958a}.
\end{proof}

We are going to apply this lemma with the choice $x= x_d y_j$ and $y = y_d x_j$ and show, in view of the lower bound \eqref{eq:BirDav:Lemma12:Count:N_j} for $N_j$, that the first alternative in the above dichotomy cannot hold. To do this, we need to adapt Lemma 13 of \cite{birch-davenport:1958a} as follows.

\begin{lemma}
  \label{lemma:BirDav:Lemma13}
  Let $j \neq l$. For any $\alpha \in \calG$ we have
  \begin{equation}
    \label{eq:BirDav:Lemma13:eq1}
    0 < \abs{x_l} y_j \ll \abs{q_l} U_l U_j u(P)
  \end{equation}
  for all integral pairs $(x_j,y_j) \in \mathfrak{D}_j(\alpha)$, $(x_l,y_l) \in \mathfrak{D}_l(\alpha)$ and also
  \begin{equation}
    \label{eq:BirDav:Lemma13:eq2}
    \big\lvert x_l y_j \frac{q_j}{q_l} - x_j y_l \big\rvert \ll \abs{q_j}  (U_l U_j)^\frac{1}{2} (T_l T_j)^{-1} P^{-2} (\log{P})^2.
  \end{equation}
\end{lemma}

\begin{proof}
   We recall that $x_i \neq 0$ for any $i=1,\ldots,d$ and that the size of $\abs{x_l}$ is of order
  \begin{equation*}
    \abs{q_l} \alpha y_l \ll \abs{x_l} \ll \abs{q_l} \alpha y_l,
  \end{equation*}
  because the approximation error $\abs{\rho_j}$ is small compared to $\abs{q_l} \alpha y_l$, see \eqref{eq:Dio:2}. Thus
  \begin{equation*}
    0 < \abs{x_l} y_j \ll \abs{q_l} \alpha y_l y_j \ll u(P) \abs{q_l} U_l U_j,
  \end{equation*}
  where  $y_j y_l \le U_j U_l$, see \eqref{eq:DefSet:G}, and $\alpha < u(P)$ was used. To prove \eqref{eq:BirDav:Lemma13:eq2}, we note first that
  \begin{equation*}
    2\alpha = \frac{1}{q_j}\frac{x_j}{y_j}+\frac{\rho_j}{q_j} = \frac{1}{q_l}\frac{x_l}{y_l}+\frac{\rho_l}{q_l}.
  \end{equation*}
  Hence after multiplying by $y_l y_j q_j$ and arranging accordingly we see that
  \begin{equation*}
    x_l y_j\frac{q_j}{q_l} - x_jy_l = y_l y_j q_j ({q_j}^{-1}{\rho_j}  -  {q_l}^{-1}{\rho_l}).
  \end{equation*}
  Consequently, as in the proof of Lemma 13 in \cite{birch-davenport:1958a}, we have
  \begin{equation*}
    \bigg|x_l y_j \frac{q_j}{q_l}- x_j y_l\bigg| \ll y_l y_j \abs{q_j} (\abs{q_j^{-1}\rho_j} + \abs{q_l^{-1} \rho_l}).
  \end{equation*}
  The inequality \eqref{eq:BirDav:eq56}, that is $\abs{q_i}^{-1} \abs{\rho_i} \ll (\log{P}) P^{-2} T_i^{-1} U_i^{-1/2}$, combined with the definition \eqref{eq:DefSet:G} of $U_j$ shows that the last term can be bounded by
  \begin{equation*}
    \ll \abs{q_j} U_l U_j \big( T_j^{-1} U_j^{-1/2} + T_l^{-1} U_l^{-1/2}\big) P^{-2} (\log{P})
  \end{equation*}
  and, in view of the relation \eqref{eq:UpperBoundForU_j} between $T_i$ and $U_i$, this is further bounded by
  \begin{equation*}
    \ll \abs{q_j} (U_l U_j)^{1/2} (T_l T_j)^{-1} P^{-2} (\log{P})^2.
  \end{equation*}
  This proves \eqref{eq:BirDav:Lemma13:eq2}.
\end{proof}

The first part of the following lemma will be essential for verifying that at least $d-3$ variables are coupled, whereas the second part will be used for both smaller dimensions and quadratic forms of signature $(r,s)$ with relatively large exponent $\beta(r,s)$ (recall that $\beta(r,s)$ was introduced in Theorem \ref{th:Intro:SchlickeweiTheorem}).

\begin{lemma}
  \label{lemma:BirDav:Lemma15}
  If $d \ge 8$ and $j \in \{4,\ldots,d-1\}$, then for any $\alpha \in \calG$ we have
  \begin{equation}
    \label{eq:Lemma15:Coulpling}
    \frac{x_j y_d}{y_j x_d} = \frac{A_j}{B_j},
  \end{equation}
  where $(x_j,y_j) \in \mathfrak{D}_j(\alpha)$, $(x_d,y_d) \in \mathfrak{D}_d(\alpha)$, and $B_j > 0, A_j \neq 0$ are coprime integers independent of $\alpha$. The same holds also for the coordinates
    \begin{enumerate}[(a)]
      \item $3 \le j \le d-1$ if $\beta \ge 2/3$ and $d \ge 7$,
      \item $2 \le j \le d-1$ if $\beta \ge 1$ and $d \ge 6$,
      \item $1 \le j \le d-1$ if $\beta \ge 2$ and $d \ge 5$.
    \end{enumerate}
\end{lemma}

 At this point the initial coupling only applies to at least $d-3$ coordinates because the size of $T_j$ can only be tamed by the first $j$ numbers $T_1, \ldots, T_j$. Due to this the second alternative in Lemma \ref{lemma:BirDav:Lemma13} can be excluded only if $j \geq 4$. The general problem is to extract from the lower bound \eqref{eq:BirDav:Lemma11}  information about the  sums $S_j$ which  only  have a small contribution to the integral under consideration.
 
 If $j \in \{1,2,3\}$, then the size of $\beta$ becomes important (note that this is the first time it arises) because the lower bound on $T_j$, stated in \eqref{eq:BoundsForT_j}, has to be used (this bound is rather weak, but the argument used there cannot be improved, see Remark \ref{remark:cannot_improve}). Since the size of the approximation error \eqref{eq:BirDav:Lemma13:eq2} crucially depends on the size of $\beta$, this is only feasible if $\beta$ is not too small.

\begin{proof}
  The general strategy here is to apply Lemma \ref{lemma:dioph_approx} to the integers $x= x_d y_j$ and $y = y_d x_j$, where $(x_j,y_j) \in \mathfrak{D}_j(\alpha)$ and $(x_d,y_d) \in \mathfrak{D}_d(\alpha)$ for some $\alpha \in \calG$. We only carry out the proof for $j \in \{4,\ldots,d-1\}$ and afterwards outline the required changes for the remaining cases (a)--(c). By Lemma \ref{lemma:BirDav:Lemma13} we have
  \begin{align*}
   \abs{x q_j/q_d - y} < \eta \quad &\text{and} \quad 0 < \abs{x} < X \quad \text{with} \\
   X \ll u(P) \abs{q_d} (U_d U_j) u(P) \quad &\text{and} \quad \eta \ll \abs{q_j} (U_d U_j)^\frac{1}{2} (T_d T_j)^{-1} P^{-2} u(P).
  \end{align*}
  According to Lemma \ref{lemma:dioph_approx} either $N \le 24 \eta X$, where $N$ denotes the number of distinct integer pairs $(x,y)$ corresponding to any $\alpha \in \calG$, or all pairs $(x,y)$ have the same ratio $y / x$, independent of $\alpha$, which gives the desired conclusion. We show that the former case is impossible, provided $C_d \gg 1$ is chosen sufficiently large: In the former case case, we have
  \begin{equation}
    \label{eq:proof:BirDav:Lemma15:eq1}
    N \leq 24 \eta X \ll \abs{q_d q_j} (U_d U_j)^\frac{3}{2} (T_d T_j)^{-1} P^{-2} u(P)^2
  \end{equation}
  and, furthermore, the values of $x_d,y_d$ are determined by the divisors of $x$ and $y$. Since there are $\ll P^\rho$ divisors (for any fixed $\rho >0$) and $x_d \neq 0$, we find
  \begin{equation*}
    N_d \ll P^\rho N.
  \end{equation*}
  Now we may use the lower bound \eqref{eq:BirDav:Lemma12:Count:N_j} from Corollary \ref{lemma:BirDav:Lemma12} together with the upper bound \eqref{eq:proof:BirDav:Lemma15:eq1} to get
  \begin{equation*}
    (\log{P})^{-2d} (T_1 \ldots T_d)^{-1} (T_d U_d^{1/2}) \ll \abs{q_1 q_j} (U_d U_j)^{3/2} (T_d T_j)^{-1} P^{-2+ \rho} u(P).
  \end{equation*}
  By \eqref{eq:UpperBoundForU_j} this can simplified as
  \begin{equation}
    \label{eq:Proof:Coupling}
    T_d^4 T_j^4 \ll q^2 P^{-2+\rho} (\log{P})^{2d+5} u(P)^2 (T_1\ldots T_d).
  \end{equation}
  Suppose that $j \in \{4, \ldots,d-1\}$ and $d \ge 8$. In this situation we have $(T_1\ldots T_d) \ll T_d^4 T_j^4$, where we used $T_1 \ll \ldots \ll T_d$ together with $T_i \ll 1$, compare \eqref{eq:BoundsForT_j}. In conjunction with \eqref{eq:Proof:Coupling} we now deduce the inequality
  \begin{equation*}
    T_d^4 T_j^4 \ll q^2 P^{-2+\rho} (\log{P})^{2d+5} u(P)^2 T_d^4 T_j^4
  \end{equation*}
  and by canceling $T_j^4$ and $T_d^4$ on both sides we further obtain
  \begin{equation}
    \label{eq:Proof:contradiction}
    1 \ll q^2 P^{-2+\rho} (\log{P})^{2d} u(P)^2 .
  \end{equation}
  Since $2\beta \ge (d+3)/(d-3)>1$, we can choose $\rho >0$ such that $2 < (2-\rho)(1+2\beta)$ and note that the right hand side of \eqref{eq:Proof:contradiction} tends to zero, compare \eqref{eq:ChoiceOfP}. Thus, after increasing $C_d \gg 1$, we find that inequality \eqref{eq:Proof:contradiction} cannot hold (note that the implicit constant depends on $d$ only). This contradiction shows that the first alternative in Lemma \ref{lemma:dioph_approx} must hold.
  
  In the other cases we should use Wigert's divisor bound, i.e.\ 
  \begin{equation}
    \label{eq:wigerts_bound}
    d(n) \ll_\eps 2^{(1+\eps) \log(n) /\log \log{n}}
  \end{equation}
  if $\eps >0$ (for a reference, see Theorem 317 in \cite{Hardy-Wright:2008}), regarding that $\abs{x},\abs{y} \ll P^3$. If $ 3 \le j \le d-1$ and $d \ge 7$, then we still find that 
  \begin{equation*}
     T_d^4 T_j^4 \ll q^2 P^{-2+\frac{6}{\log \log P}} (\log{P})^{2d+5} u(P)^2 T_d^4 T_j^3.
  \end{equation*}
  Canceling $T_d^4 T_j^3$ and using $T_j \gg q^{-1/3} u(P)^{-2/3}$, compare \eqref{eq:BoundsForT_j}, gives
  \begin{equation*}
    1 \ll q^{7/3} P^{-2+\frac{6}{\log \log P}} (\log{P})^{2d+5} u(P)^{8/3}.
  \end{equation*}
  To deduce a contradiction again, we need at least $1+2 \beta \geq 7/3$ (here the precise definition of $P$ in terms of $H$ is used, compare \eqref{eq:ChoiceOfP}). The remaining cases can be proved similarly: If $2 \le j \le d-1$ and $d \ge 6$, then we obtain the inequality
  \begin{equation*}
    T_j^2 \ll q^2 P^{-2+\frac{6}{\log \log P}} (\log{P})^{2d+5} u(P)^2.
  \end{equation*}
  By using $T_j \gg q^{-1/2} u(P)^{-1}$ we see that at least $1+ 2\beta \ge 3$ is required. In the last case, i.e.\ $1 \le j \le d-1$ and $d \ge 5$, we need $1+ 2\beta \ge 5$, since we only know that 
  \begin{equation*}
    T_j^3 \ll q^2 P^{-2+\frac{6}{\log \log P}} (\log{P})^{2d+5} u(P)^2 
  \end{equation*}
  and $T_j \gg q^{-1} u(P)^{-2}$.
\end{proof}

The above lemma allows us to obtain a factorization of $x_j$ and $y_j$ as formulated in the Definition \ref{def:coupling} of the notion of `coupling'.

\begin{lemma}
  \label{lemma:BirDav:Lemma16}
  Let $I \subset \{1,\ldots,d-1\}$ be some set of indices. Assume that for any $\alpha \in \mathcal{G}$ and all $j \in I$ the integral pairs $(x_j,y_j) \in \mathfrak{D}_j(\alpha), (x_d,y_d) \in \mathfrak{D}_d(\alpha)$ can be factorized as in \eqref{eq:Lemma15:Coulpling}, where $A_j,B_j$ are coprime integers which are independent of $\alpha$ (but may depend on $I$) and $B_j > 0$, $A_j \neq 0$. Then all coordinates from the set $I \cup \{d\}$ are coupled on $\mathcal{G}$ (in the sense of Definition \ref{def:coupling}) with corresponding common multiple $L \in \N$ satisfying
  \begin{equation}
    \label{eq:Lemma16:eq3}
    0 < L \ll H^{10d},
  \end{equation}
  where $H$ is as in \eqref{eq:ChoiceOfP}.
\end{lemma}

\begin{proof}
  For any $j \in I$ we can rewrite equation \eqref{eq:Lemma15:Coulpling} as 
  \begin{equation*}
    \frac{x_j}{y_j} = \frac{x_d}{y_d} \frac{B_j}{A_j},
  \end{equation*}
  where $(x_d,y_d) = (x_j,y_j) = (A_j,B_j)=1$. This reveals the factorization
  \begin{equation*}
    x_j = \sgn(x_j) \frac{\abs{x_d}}{(x_d,A_j)} \frac{B_j}{(y_d,B_j)} \quad \text{and} \quad y_j = \frac{y_d}{(y_d,B_j)} \frac{\abs{A_j}}{(x_d,A_j)}
  \end{equation*}
  and this factorization allows us to define $x$ and $y$ by
  \begin{equation*}
    x = \frac{\abs{x_d}}{(x_d, \prod_{i \in I} A_i)} \quad \text{and} \quad y = \frac{y_d}{(y_d, \prod_{i \in I} B_i)}.
  \end{equation*}
  Then $x$ and $y$ are non-zero integers and we can further define
  \begin{equation*}
    x_d' := \frac{x_d}{x} = \sgn(x_d) (x_d, \mprod_{i \in I} A_i) \quad \text{and} \quad y_d' := \frac{y_d}{y} = (y_d, \mprod_{i \in I} B_i)
  \end{equation*}
  and also
  \begin{equation*}
    x_j' := \frac{x_j}{x} = \sgn(x_j) \frac{\abs{x_d} B_j}{(x_d,A_j)(y_d,B_j)} \frac{(x_d, \mprod_{i \in I} A_i)}{\abs{x_d}}
    = \sgn(x_j) \frac{B_j}{(y_d,B_j)} \frac{(x_d, \mprod_{i \in I} A_i)}{(x_d,A_j)}
  \end{equation*}
  and
  \begin{equation*}
    y_j' := \frac{y_j}{y} = \frac{y_d \abs{A_j}}{(y_d,B_j)(x_d,A_j)} \frac{(y_d, \mprod_{i \in I} B_i)}{y_d} = \frac{\abs{A_j}}{(y_d,A_j)} \frac{(y_d, \mprod_{i \in I} B_i)}{(y_d,B_j)}.
  \end{equation*}
  Note that both are non-zero integral numbers and that $x_j'$ and $y_j'$ are divisors of
  \begin{equation*}
    L:= \mprod_{i \in I} \abs{A_i B_i}.
  \end{equation*}
  It remains to find an upper bound for $K$. By Lemma \ref{lemma:BirDav:Lemma13} we have
  \begin{equation*}
    \abs{A_j B_j} \leq \abs{x_d} y_j \abs{x_j} y_d \ll u(P)^2 \abs{q_1} \abs{q_j} U_d^2 U_j^2
  \end{equation*}
  and thus
  \begin{equation*}
    L \ll u(P)^{2(d-1)} q^{2(d-1)} U_d^{2(d-1)} \mprod_{i \in I} U_i^2 \ll u(P)^{2(d-1)} (\log{P})^{4(d-1)} q^{2(d-1)} T_d^{-4(d-1)} \mprod_{i \in I} T_i^{-4},
  \end{equation*}
  where we used \eqref{eq:UpperBoundForU_j}. In view of \eqref{eq:BoundsForT_j} this is bounded by
  \begin{equation*}
    \ll u(P)^{2(d-1)} \log(P)^{4(d-1)} q^{2(d-1)} (u(P)^2 q)^{\frac{4(d-1)}{(d-4)}+ 4\sum_{i=1}^{d-1} \kappa(i)} \ll u(P)^{16d+8 \log(d)} q^{8d + 4\log(d)}.
  \end{equation*}
  Using the definition of $H$ (see \eqref{eq:ChoiceOfP}) together with $\beta \geq 1/2$ yields that the last inequality chain is at most $\ll H^{10d}$.
\end{proof}

 Combining Lemmata \ref{lemma:BirDav:Lemma15} and \ref{lemma:BirDav:Lemma16} shows that in each of the cases of Lemma \ref{lemma:BirDav:Lemma15} all indices $j$ under consideration are coupled on $\mathcal{G}$ with a common multiple $L \in \N$ (as in the Definition \ref{def:coupling}) being bounded as in \eqref{eq:Lemma16:eq3}. Additionally, taking into account the definition of $\beta(r,s)$ for a given dimension $d$ and given signature $(r,s)$, we conclude the following

\begin{corollary}
  \label{cor:first_coupling}
  Under Assumption \ref{assumption:Theorem1} the functions $S_4, \ldots, S_d$ are always coupled on $\mathcal{G}$. Assuming additionally the following conditions imply that $S_{k+1}, \ldots, S_{d}$ are coupled on $\mathcal{G}$.
  \begin{enumerate}
    \item $k=0$ if $d \in \{5,6\}$ or $r \ge 4 s$,
    \item $k=1$ if $5 \le d \le 10$ or $r \ge 2s$ and $d \ge 11$,
    \item $k=2$ if $5 \le d \le 22$ or $r \ge 4s/3$ and $d \ge 23$.
  \end{enumerate}
\end{corollary}

\section{The Rational Case: Schlickewei's Bound on Small Zeros}
\label{sec:rational-case}
\noindent In this section we shall state for the reader's convenience the main results of Schlickewei's work \cite{Schlickewei:1985} on small zeros of integral quadratic forms and afterwards deduce Theorem \ref{th:Intro:SchlickeweiTheorem} from Schlickewei's results. We also note that our proof of Theorem \ref{MainTheorem} and the exponent in \eqref{eq:Theorem1} depend essentially on the results presented here.

\begin{theorem}
  \label{theorem:Schlickewei:1}
  Let $F$ be a non-trivial quadratic form in $d$ variables with integral coefficients and let $G$ be a positive definite quadratic form with real-valued coefficients. Furthermore, let $d_0$ be maximal such that $F$ vanishes on a rational subspace of dimension $d_0$. Then there exist integral points $\mathfrak{M}_1,\ldots,\mathfrak{M}_{d_0}$, which are linearly independent over $\Q$, such that $F$ vanishes on the spanned $\mathbb{Q}$-subspace and
  \begin{equation}
    \label{eq:Schlickewei:eq1}
    0 < G(\mathfrak{M}_1) \cdot \ldots \cdot G(\mathfrak{M}_{d_0}) \ll_d \mathrm{trace}((FG^{-1})^2)^\frac{d-d_0}{2} \det{G},
  \end{equation}
  where the constant in $\ll$ depends on the dimension $d$ only.
\end{theorem}

\begin{remark-non}
 For the latter application of this bound the explicit dependence on the determinant is crucial. One of the reasons for this is that the lower bound \eqref{eq:BirDav:Lemma12:Count:N_j} on $N_j$ can also be written in terms of the determinant of $\operatorname{diag}(T_1,\ldots,T_d)$.
\end{remark-non}

Theorem \ref{theorem:Schlickewei:1} is \textsl{Satz 2} in \cite{Schlickewei:1985} and the proof relies on an application of Minkowski's second theorem on successive minima. Moreover, by using an induction argument combined with Meyer's theorem \cite{Meyer:1884}, Schlickewei found the following connection between the dimension of a maximal rational isotropic subspace and the signature.

\begin{proposition}
  \label{theorem:Schlickewei:2}
  Let $F$ be a quadratic form in $d$ variables with integral coefficients and signature $(r,s,t)$. Suppose that $r \geq s$ and $r+s \geq 5$. The dimension $d_0$ of a maximal rational isotropic subspace is at least
  \begin{equation}
    \label{eq:Schlickewei:35}
    d_0 \geq \begin{cases}
      s+t   & \text{if} \ r \geq s+3                \\
      s+t-1 & \text{if} \ r=s+2 \ \text{or} \ r=s+1 \\
      s+t-2 & \text{if} \ r=s
    \end{cases}.
  \end{equation}
\end{proposition}

Note that the quadratic form $F$ is allowed to be degenerate and then the triple $(r,s,t)$ expresses the number $r$ of positive, $s$ of negative and $t$ of zero entries in its reduced form.

\begin{proof}
  See in \cite{Schlickewei:1985}, \textsl{Hilfssatz} in Section 4.
\end{proof}

Now, we can argue like Schlickewei in \cite{Schlickewei:1985}, see \textit{Folgerung 3}, to deduce

\begin{corollary}
  \label{lemma:BirDav:Lemma17}
  For any non-zero integers $f_1,\ldots,f_d$, of which $r \geq 1$ are positive and $s \geq 1$ negative with $r \ge s$, $d=r+s \ge 5$, there exist integers $m_1,\ldots,m_d$, not all zero, such that
  \begin{equation}
    \label{eq:consequence-schlickewei}
    \begin{gathered}
    f_1 m_1^2 + \ldots + f_d m_d^2 =0, \\
    0 < \abs{f_1} m_1^2 + \ldots + \abs{f_d} m_d^2 \ll_d \abs{f_1 \ldots f_d}^\frac{2\beta+1}{d},
    \end{gathered}
  \end{equation}
  where $\beta$ is defined as in \eqref{eq:Intro:Beta} and the implicit constant depends on the dimension $d$ only.
\end{corollary}

\begin{proof}
  We apply Theorem \ref{theorem:Schlickewei:1} to the forms $F(m) = \sum_{j=1}^d f_j m_j^2$ and $G(m) = \sum_{j=1}^d \abs{f_j} m_j^2$ to get isotropic integral points $\mathfrak{M}_1,\ldots,\mathfrak{M}_{d_0}$ satisfying \eqref{eq:Schlickewei:eq1}. Let $\mathfrak{M}_i = (m_1,\ldots,m_d)$ be a point with minimal weight, i.e.\ $g(\mathfrak{M}_i) = \min_{j \in \{1,\ldots,d_0\}} g(\mathfrak{M}_j)$. This lattice point satisfies
  \begin{equation*}
    \abs{f_1}m_1^2 + \ldots + \abs{f_d} m_d^2 \ll \abs{f_1 \ldots f_d}^{1/d_0}.
  \end{equation*}
  If $r \geq s+3$ or $r=s+1$, we can use the lower bound \eqref{eq:Schlickewei:35} for $d_0$. But if $r=s+2$ or $r=s$, then set one variable $x_i$, such that $\abs{f_i}$ is maximal, to zero: It follows that $F$ has signature $(r,s-1)$ or $(r-1,s)$. The previous argument (applying the lower bound \eqref{eq:Schlickewei:35} again) together with the estimate
  \begin{equation*}
    \textstyle
    \prod_{j \ne i} \abs{f_j} \le \abs{f_1 \ldots f_d}^{(d-1)/d}
  \end{equation*}
  implies the claimed bound \eqref{eq:consequence-schlickewei} in these cases as well.
\end{proof}

\section{Iteration of the Coupling Argument}
\label{sec:iteration}

\noindent In Section \ref{subec:coupling} we showed that the functions $S_{k+1}, \dots, S_{d}$ are coupled on $\calG$ for some $k \in \{0,1,2,3\}$ depending on the exponent $\beta(r,s)$ introduced in Theorem \ref{th:Intro:SchlickeweiTheorem}. Namely, the integer pairs $(x_i,y_i) \in \mathfrak{D}_i(\alpha)$ corresponding to $q_i \alpha$ and any $i = k{+}1,\ldots,d$ are of the form
\begin{equation}
  \label{def:factorization}
  x_{i}= x \4 x_i' \quad \quad \text{and} \quad \quad y_{i}= y \4 y_i',
\end{equation}
with $x>0$, $y>0$, $x_{i}$ and $x_i'$ have the same sign, $x_i ' \mid L$ and $y_i' \mid L$, where $L$ is independent of $\alpha \in \calG$ and $L \ll H^{10d}$. In this section we shall utilize this observation in combination with Schlickewei's bound on small zeros in order to count the number of distinct pairs $(x,y)$. For this purpose, we introduce the set $\mathfrak{C}_k(\alpha)$ of all pairs $(x,y)$ corresponding to some fixed $\alpha \in \calG$ and here we shall always assume that $x_i$ and $y_i$ are factorized as in \eqref{def:factorization} without mentioning this explicitly.

\begin{lemma}
  \label{lemma:bound_for_xy}
  Suppose that the exponential sums $S_{k+1},\ldots,S_{d}$ are coupled on $\mathcal{G}$, where $k \in \{0,1,2,3\}$, and that the quadratic form
  \begin{equation}
    \label{def:Q_j}
    Q_k[m] := q_{k+1} m_{k+1}^2 + \ldots + q_{d} m_{d}^2
  \end{equation}
  is indefinite of signature $(r',s')$ with $d-k \ge 5$. Then, under Assumption \ref{assumption:Theorem1}, the integer pairs $(x, y) \in \mathfrak{C}_k(\alpha)$, corresponding to the factorization \eqref{def:factorization} and any $\alpha\in \calG$, satisfy
  \begin{equation}
    \label{eq:bound_for_xy:eq0}
    x^{2 \beta_k} y^{2 \beta_k +2} \ll  q^{2\beta_k+1} P^{-2} (\log{P}) u(P)^{2\beta_k} (U_{k+1} \ldots U_{d})^{\frac{4\beta_k+2}{d-k}} \big(\max_{i=k+1,\ldots,d} T_i^{-1} U_i^{-\frac{1}{2}}\big),
  \end{equation}
  where $\beta_k = \beta(r',s')$ denotes the exponent (as defined in \eqref{eq:Intro:Beta} of Theorem \ref{th:Intro:SchlickeweiTheorem}) corresponding to the signature $(r',s')$ of $Q_k$ and $u(P)$ is chosen as in \eqref{eq:choice_u}.
\end{lemma}

This lemma will be used subsequently to establish improved mean value estimates and, as a consequence, improved lower bounds for the size of the parameters $T_1,\ldots,T_k$.

\begin{proof}
  Due to the Diophantine approximation introduced in \eqref{eq:Dio:1}, we have for any fixed $\alpha \in \calG$ and any integers $m_{k+1},\ldots,m_{d} \in \Z$
  \begin{equation*}
    \alpha(q_{k+1} m_{k+1}^2 + \ldots + q_{d} m_{d}^2) = \frac{x}{y} \bigg( \frac{x_{k+1}'}{y_{k+1}'} m_{k+1}^2 + \ldots + \frac{x_{d}'}{y_{d}'} m_{d}^2 \bigg) + (\rho_{k+1} m_{k+1}^2 + \ldots + \rho_{d} m_{d}^2).
  \end{equation*}
  Here we change variables to $m_i = y_i' n_i$ for any $i=k+1,\dots,d$ and get
  \begin{equation}
    \label{eq:Proof:Lemma18:eq1}
    \alpha(q_{k+1} m_{k+1}^2 + \ldots + q_{d} m_{d}^2)  = \frac{x}{y} \sum_{j=k+1}^d x_{j}' y_{j}' n_{j}^2 + \sum_{j=k+1}^d \rho_{j} y_{j}'^2 n_{j}^2.
  \end{equation}
  Observe that the first term on the right hand side, neglecting the factor $x/y$, is an integral quadratic form whose signature $(r',s')$ coincides with that of $Q_k$, since the signs of $x_{k+1}'y_{k+1}',\ldots,x_{d}'y_{d}'$ are exactly equal to those of $x_{k+1}/y_{k+1},\ldots,x_{d}/y_{d}$ and these have the same signs as $q_{k+1},\ldots,q_{d}$. Hence, it follows from Corollary \ref{lemma:BirDav:Lemma17} that there exist integers $n_{k+1},\ldots,n_{d}$, not all zero, such that
  \begin{equation*}
    x_{k+1}'y_{k+1}' n_{k+1}^2 + \ldots + x_{d}' y_{d}' n_{d}^2 = 0
  \end{equation*}
  and
  \begin{equation}
    \label{eq:Proof:Lemma18:eq2}
    \abs{x_{k+1}'y_{k+1}'} n_{k+1}^2 + \ldots + \abs{x_{d}' y_{d}'} n_{d}^2 \ll_d \abs{x_{k+1}' y_{k+1}' \ldots x_{d}' y_{d}'}^{(2\beta_k+1)/(d-k)}.
  \end{equation}
  For the corresponding $m_{k+1},\ldots,m_{d}$ the first part of the right hand side in \eqref{eq:Proof:Lemma18:eq1} vanishes. Thus, we find
  \begin{equation*}
    \lvert q_{k+1} m_{k+1}^2+ \ldots + q_{d} m_{d}^2 \rvert \ll \alpha^{-1} (\abs{\rho_{k+1}} y_{k+1}'^2 n_{k+1}^2 + \ldots + \abs{\rho_{d}} y_{d}'^2 n_{d}^2)
  \end{equation*}
  and from $\alpha \lvert q_i \rvert \ll \lvert x_i \rvert y_i^{-1}$, \eqref{eq:Proof:Lemma18:eq2} and $\lvert x_{k+1}' y_{k+1}' \rvert \ll (xy)^{-1} \lvert q_i \rvert \alpha^{-1} y_j^2$ we deduce that
  \begin{equation}
  \begin{aligned}
\label{eq:Proof:Lemma18:eq3}
    \abs{q_{k+1}} m_{k+1}^2 + \ldots + \abs{q_{d}} m_{d}^2 & \ll \alpha^{-1} x y^{-1} \abs{x_{k+1}' y_{k+1}' \ldots x_{d}' y_{d}'}^{(2\beta_k+1)/(d-k)}                                \\
    & \ll \alpha^{2 \beta_k } x^{-2\beta_k} y^{-2\beta_k -2} \abs{q_{k+1} y_{k+1}^2 \ldots q_{d} y_{d}^2}^{(2\beta_k+1)/(d-k)} \\
    & \ll \alpha^{2 \beta_k} x^{-2\beta_k} y^{-2\beta_k -2} q^{2\beta_k+1} (U_{k+1} \ldots U_{d})^{(4\beta_k+2)/(d-k)},
  \end{aligned}
  \end{equation}
  where $y_i \le U_i$ was used in the last step. Now we shall apply the Assumption \ref{assumption:Theorem1}, made at the beginning: Since $Q_k$ is a restriction of $Q$, i.e.\ $Q_k[m] = Q[(0,\ldots,0,m_{k+1},\ldots,m_{d})]$, we have either
  \begin{equation}
    \label{proof:firstcase}
    4d^3P^2 <\abs{q_{k+1}} m_{k+1}^2 + \ldots + \abs{q_{d}} m_{d}^2
  \end{equation}
  or
  \begin{equation}
    \label{proof:secondcase}
    1 \leq \abs{q_{k+1} m_{k+1}^2+ \ldots + q_{d} m_{d}^2} \le \alpha^{-1} (\abs{\rho_{k+1}} y_{k+1}'^2 n_{k+1}^2 + \ldots + \abs{\rho_{d}} y_{d}'^2 n_{d}^2).
  \end{equation}
  In the first case we may combine \eqref{proof:firstcase} together with \eqref{eq:Proof:Lemma18:eq3} to get
  \begin{equation*}
    P^2 \ll \alpha^{2 \beta_k} x^{-2\beta_k} y^{-2\beta_k -2} q^{2\beta_k+1} (U_{k+1} \ldots U_{d})^{(4\beta_k+2)/(d-k)}
  \end{equation*}
  and in view of \eqref{eq:UpperBoundForU_j}, that is $T_i^{-1} U_i^{-1/2} \gg  \log{P}$, together with $\alpha < u(P)$ we conclude already that inequality \eqref{eq:bound_for_xy:eq0} holds. In the second case \eqref{proof:secondcase} holds and here we use \eqref{eq:BirDav:eq56}, i.e.\ $\abs{\rho_i} \ll \abs{q_i} P^{-2} (\log{P}) T_i^{-1} U_i^{-1/2}$, to obtain
  \begin{equation*}
    1 \ll  \alpha^{-1} \sum_{j=k+1}^d \abs{\rho_{j}} y_{j}'^2 n_{j}^2 \ll \alpha^{-1} P^{-2} (\log{P}) \big( \max_{i=k+1,\ldots,d} T_i^{-1} U_i^{-1/2} \big) \Big( \sum_{i=k+1}^d \abs{q_{i}} m_{i}^2 \Big),
  \end{equation*}
 which implies together with  \eqref{eq:Proof:Lemma18:eq3}
  \begin{equation*}
    1 \ll \alpha^{2\beta_k-1} x^{-2\beta_k} y^{-2\beta_k-2} q^{2\beta_k+1} P^{-2} (\log{P}) \big(\max_{i=k+1,\ldots,d}  T_i^{-1} U_i^{-1/2}\big) (U_{k+1} \ldots U_{d})^{(4\beta_k+2)/(d-k)}.
  \end{equation*}
  Finally, taking into account that $2\beta_k \ge 1$ and $\alpha < u(P)$ proves inequality \eqref{eq:bound_for_xy:eq0}.
\end{proof}

All pairs $(x,y) \in \mathfrak{C}_k := \{(x,y) \in \Z^2 _{\text{prim}} : \4 (x,y) \in \mathfrak{C}_k(\alpha) \text{ for some } \alpha \in \calG\}$ lie in a bounded set de\-ter\-mined by condition \eqref{eq:bound_for_xy:eq0}. Hence, we can bound the number $\# \mathfrak{C}_k$ of all these pairs as follows.

\begin{corollary}
  \label{cor:bound_fundamental}
  In the situation of Lemma \ref{lemma:bound_for_xy}, we have
  \begin{equation}
     \label{cor:bound_fundamental:eq}
     \# \mathfrak{C}_k \ll q^{1+\frac{1}{2\beta_k}} P^{-\frac{1}{\beta_k}} (\log{P})^{\frac{1}{2\beta_k}} u(P) ( U_{k+1} \ldots U_d)^{\frac{4\beta_k+2}{2\beta_k(d-k)}} \big(\max_{i=k+1,\ldots,d} T_i^{-1} U_i^{-\frac{1}{2}}\big)^{\frac{1}{2\beta_k}}.
  \end{equation}
\end{corollary}

\begin{proof}
 First note that the expression on the right hand side of \eqref{eq:bound_for_xy:eq0} must be $\gg 1$, since $\calG$ is not empty. Thus, we can apply Dirichlet's hyperbola method to see that the number $N$ of distinct solutions $(x,y)$ of
  \begin{equation*} 
    x^{2\beta_k} y^{2\beta_k+2} \ll Z
  \end{equation*}
  is $\ll Z^\frac{1}{2\beta_k}$. This already concludes the proof.
\end{proof}

We are in position to establish improved mean value estimates (conditionally under Assumption \ref{assumption:Theorem1}) by controlling the sum over all $(x,y) \in \mathfrak{C}_k$ with the help of Corollary  \ref{cor:bound_fundamental}.


\begin{lemma}
  \label{averagewithresonance}
  Suppose that $d \ge 5+k$ and $k \in \{0,1,2,3\}$. Then for any $\rho >0$, in the situation of Lemma \ref{lemma:bound_for_xy}, we have
  \begin{equation}
    \label{eq:averagewithresonance}
    \int_{\calG} \abs{S_{k+1}(\alpha) \ldots S_{d}(\alpha) K(\alpha)} \, \dif \alpha \ll_{\rho} P^\rho \frac{P^{d-k-2}}{\abs{q_{k+1} \ldots q_{d}}^{\frac{1}{2}}} \frac{q^{1+\frac{1}{2\beta_k}}}{P^{\frac{1}{\beta_k}}}.
  \end{equation}
\end{lemma}

\begin{proof}
  We shall decompose the integration domain $\calG$ according to the covering induced by the factorization from \eqref{def:factorization}, which holds since $S_{k+1},\ldots,S_{d}$ are coupled on $\mathcal{G}$: For fixed $(x,y) \in \mathfrak{C}_k$ we define
  \begin{align*}
    \mathfrak{H}_i(x,y) := \{ (x_i',y_i') \in \Z^2 _{\text{prim}} \, : \, x_i = x x_i'& \text{ and } y_i = y y_i' \text{ as in } \eqref{def:factorization} \\
     &\text{ with } (x_i,y_i) \in \mathfrak{D}_i(\alpha) \text{ for some } \alpha \in \calG \}
  \end{align*}
  and
  \begin{equation*}
    \calJ_i(x_i,y_i) := \{ \alpha \in \calG \, : \, \abs{\alpha q_i y_i - x_i} < \abs{q_i}^{1/2} (8dP)^{-1} \}
  \end{equation*}
  in order to obtain the decomposition
  \begin{equation*}
    \operatorname{LHS} \eqref{eq:averagewithresonance} \leq \sum_{(x,y) \in \mathfrak{C}_k} \sum_{(x_{k+1}',y_{k+1}') \in \mathfrak{H}_{k+1}(x,y) } \! \ldots \! \sum_{(x_{d}',y_{d}') \in \mathfrak{H}_{d}(x,y) } I(x_{k+1},y_{k+1},\ldots,x_{d},y_{d}),
  \end{equation*}
  where
  \begin{equation*}
    I(x_{k+1},y_{k+1},\ldots,x_{d},y_{d}) := \int_{ \bigcap_{i=k+1}^{d} \calJ_i(x_i,y_i) } \lvert S_{k+1}(\alpha) \ldots S_{d}(\alpha) K(\alpha) \rvert \, \dif \alpha.
  \end{equation*}
  Using the bound $\lvert S_i(\alpha) \rvert \le \lvert q_i \rvert^{-1/2} T_i P$, compare the definition \eqref{eq:DefSet:G} of the set $\calG$, yields
  \begin{equation*}
    I(x_{k+1},y_{k+1},\ldots,x_{d},y_{d}) \le \frac{P^{d-k}(\log{P})}{\lvert q_{k+1} \ldots q_{d}\rvert^{1/2}}  (T_{k+1} \ldots T_{d}) \operatorname{mes}( \mcap_{i=k+1}^{d} \calJ_i(x_i,y_i))
  \end{equation*}
  and, since the measure of the set $\calJ_i(x_i,y_i)$ is at most $\ll P^{-2} (\log{P}) T_i^{-1} U_i^{-1/2}$, H\"{o}lder's inequality implies
  \begin{equation*}
   I(x_{k+1},y_{k+1},\ldots,x_{d},y_{d}) \ll \frac{P^{d-k-2} (\log{P})}{\lvert q_{k+1} \ldots q_{d} \rvert^{1/2}}  (T_{k+1} \ldots T_{d}) \mprod_{i=k+1}^{d} (T_i^{-1} U_i^{-1/2})^{\frac{1}{d-k}}.
  \end{equation*}
  Returning to the initial decomposition of the integral, we note that $\# \mathfrak{H}_i(x,y) \ll P^\rho$, because $x_i',y_i'$ are divisors of $L \ll H^{10d}$ and there are at most $\ll P^\rho$ divisors. Thus, taking all together we find
  \begin{equation*}
     \operatorname{LHS} \eqref{eq:averagewithresonance} \ll P^\rho \frac{P^{d-k-2} (\log{P})}{\abs{q_{k+1} \ldots q_{d}}^{1/2}}  (T_{k+1} \ldots T_{d}) \big(\mprod_{i=k+1}^{d} (T_i^{-1} U_i^{-1/2})^{\frac{1}{d-k}} \big) \# \mathfrak{C}_k.
  \end{equation*}
  Next we insert the bound \eqref{cor:bound_fundamental:eq}, established in Corollary \ref{cor:bound_fundamental}, and conclude that the last equation is bounded by
  \begin{equation*}
  \ll P^\rho \frac{P^{d-k-2}}{\abs{q_{k+1} \ldots q_{d}}^{1/2}} \frac{q^{1+\frac{1}{2\beta_k}}}{P^{\frac{1}{2\beta_k}}} (\log{P})^2 u(P)  \big(\max_{i=k+1,\ldots,d} T_i^{-1} U_i^{-\frac{1}{2}}\big)^{\frac{1}{2\beta_k}} \prod_{i=k+1}^{d} (T_i U_i^{1/2})^{1-\frac{1}{d-k}},
  \end{equation*}
  where we used that $\frac{4\beta_k+2}{2\beta_k (d-k)} \le \frac{1}{2}$ holds provided that $d \ge 5+k$. The claim follows now from the fact that $\frac{1}{2\beta_k}+\frac{1}{d-k}-1\leq -\frac{6}{d-k+3}+\frac{1}{d-k}\leq 0$ and \eqref{eq:UpperBoundForU_j}.
\end{proof}

\begin{corollary}
  \label{cor:averagewithresonance}
  In the situation of Lemma \ref{averagewithresonance}, we have for any $\rho>0$
  \begin{equation}
    \label{cor:averagewithresonance:eq}
    T_1 \dots T_k  \gg  P^{-\rho} P^{\frac{1}{\beta_k}} q^{-1-\frac{1}{ 2\beta_k}}.
  \end{equation}
\end{corollary}

\begin{proof}
  We recall the lower bound
\begin{equation}
    \label{proof:averagewithresonance:eq1}
    \int_{\calG} \abs{S_1(\alpha) \ldots S_d (\alpha) K(\alpha)} \, \dif \alpha \gg \abs{Q}^{-\frac{1}{2}} P^{d-2} (\log{P})^{-2d}
  \end{equation}
  obtained in Lemma \ref{lemma:BirDav:Lemma11} under Assumption \ref{assumption:Theorem1}. Combining \eqref{proof:averagewithresonance:eq1} together with
  \begin{equation*}
    \abs{S_1(\alpha)\dots S_k(\alpha)} \leq  \abs{q_{1} \ldots q_k}^{-1/2} \,P^k \, (T_{1} \ldots T_{k}),
  \end{equation*}
  where we used the localization introduced in \eqref{eq:DefSet:G}, and the mean value estimate derived in Lemma \ref{averagewithresonance} shows that
  \begin{equation*}
    \abs{Q}^{-\frac{1}{2}} P^{d-2} (\log{P})^{-2d} \ll  P^{\rho/2} \abs{Q}^{-\frac{1}{2}} P^{d-2} q^{1+\frac{1}{2\beta_k}}{P^{-\frac{1}{\beta_k}}} (T_{1} \ldots T_{k}). \qedhere
  \end{equation*}
\end{proof}

\subsection{Reducing Variables and Corresponding Signatures}
\label{subsec:couple_all}
\noindent Now we are in position to prove that the remaining coordinates are coupled as well: Beginning with $S_{3}$, we will repeat the basic strategy used in the proof of Lemma \ref{lemma:BirDav:Lemma15}, but we additionally utilize the bound \eqref{cor:averagewithresonance:eq}. Compared to the earlier arguments, we need also to consider ratios between $\beta$ and $\beta_k$ with care, since simple bounds on $\beta_k$ (resp.\ on $\beta$) are not sufficient to deduce a contradiction. This step has been moved to \hyperref[appendix-a]{Appendix A}, where we address the problem to specify the possible values of $\beta_k$ depending on the signature $(r,s)$ of $Q$.

\begin{lemma}
\label{firstcoupling}
  Let $d \ge 8$ and assume that the signature of $Q$ is not of the form $(d-1,1)$, $(d-2,2)$ or $(d-3,3)$. Then, under Assumption \ref{assumption:Theorem1}, $S_{3},\ldots,S_{d}$ can be coupled on $\calG$.
\end{lemma}

\begin{proof}
 According to Corollary \ref{cor:first_coupling} we may assume that $S_{4},\ldots, S_{d}$ are coupled on $\mathcal{G}$. Applying Lemma \ref{lemma:dioph_approx} to the integers $x = x_d y_{3}$ and $y = y_d x_{3}$ with $(x_d,y_d) \in \mathfrak{D}_d(\alpha)$ and $(x_{3},y_{3}) \in \mathfrak{D}_{3}(\alpha)$ and assuming that the first alternative of Lemma \ref{lemma:dioph_approx} holds, yields (as in the proof of Lemma \ref{lemma:BirDav:Lemma15}) inequality \eqref{eq:Proof:Coupling}, that is
 \begin{equation*}
    T_d^4 T_3^4 \ll q^2 P^{-2+\rho} (\log{P})^{2d+5} u(P)^2 (T_1\ldots T_d) \ll q^2 P^{-2+\rho} (\log{P})^{2d+5} u(P)^2 (T_d^4 T_3^3),
 \end{equation*}
 where in the last step we used that $T_1 \ll \ldots \ll T_d$ and $T_i \ll 1$. Now we can cancel $T_d^4 T_{3}^3$ and use Corollary \ref{cor:averagewithresonance} with $k=3$ (note that the assumption on the signature guarantees that the quadratic form \eqref{def:Q_j} of Lemma \ref{lemma:bound_for_xy} is indefinite) to obtain
 \begin{equation*}
    P^{\frac{1}{3\beta_3}-\frac{\rho}{3}} q^{-\frac{1}{3}-\frac{1}{6 \beta_3}} \ll T_{3} \ll q^2 P^{-2+\rho} (\log{P})^{2d+5}.
 \end{equation*}
 Rearranging the last inequality and using that $q \ll P^{\frac{2}{1+2\beta}}$ gives
 \begin{equation}
    \label{eq:proof:couple_d-2}
    1 \ll P^{2\rho} (\log P)^{2d+5}u(P)^2 P^{\mathfrak p_3(d)},
 \end{equation}
 where
 \begin{equation*}
   \mathfrak{p}_3(d) := \mfrac{2}{(1+2\beta)}\big(\mfrac{7}{3}+\mfrac{1}{6 \beta_3}\big)-\big(2+\mfrac{1}{3\beta_3}\big).
 \end{equation*}
 Considering all cases in Table \ref{table:third} of the \hyperref[appendix-a]{Appendix A} we see that $\mathfrak p_3(d) <0$ and thus inequality \eqref{eq:proof:couple_d-2} cannot hold if we increase $C_d >1$ and choose $\rho>0$ small enough. (Note that Corollary \ref{cor:averagewithresonance} holds for any $\rho>0$.) To sum up, we showed that the second alternative in Lemma \ref{lemma:dioph_approx} holds, i.e.\ there exists a factorization
 \begin{equation*}
    \frac{x_3 y_d}{y_3 x_d} = \frac{A_3}{B_3}
 \end{equation*}
 for all $(x_3,y_3) \in \mathfrak{D}_3(\alpha)$, $(x_d,y_d) \in \mathfrak{D}_d(\alpha)$ and any $\alpha \in \mathcal{G}$, where $A_3,B_3$ are coprime integers which are independent of $\alpha$ and $B_3 > 0, A_3 \neq 0$. Finally, to conclude that the functions $S_3,\ldots,S_d$ are coupled on $\mathcal{G}$, note that the assumptions of Lemma \ref{lemma:BirDav:Lemma16} are now satisfied for the choice $I= \{3,\ldots,d-1\}$. (The reader may note that in each iteration step the factorization \eqref{eq:Lemma15:Coulpling} may change if further coordinates are coupled.)
\end{proof}

To proceed we need to recall some consequences of Corollary \ref{cor:first_coupling}: If $d \in \{5,6\}$, then part (i) implies that all exponential sums $S_1,\ldots,S_d$ are coupled. If $d \geq 7$, then part (iii) implies that $S_{3},\ldots S_{d}$ are coupled if $5 \leq d \leq 22$ or if $Q$ has signature $(d-1,1)$, $(d-2,2)$ or $(d-3,3)$, since in this cases $r \geq 4s/3$ is satisfied for $d \geq 23$. Hence, in view of the previous lemma, we conclude that $S_{3},\ldots,S_{d}$ are always coupled on $\mathcal{G}$.

\begin{lemma}
  Let $d \ge 7$ and assume that the signature of $Q$ is not of the form $(d-1,1)$ or $(d-2,2)$. Then, under Assumption \ref{assumption:Theorem1}, $S_{2},\ldots,S_{d}$ can be coupled on $\calG$.
\end{lemma}

\begin{proof}
 Based on inequality \eqref{eq:Proof:Coupling} we find again that $T_{2}^2 \ll q^2 P^{-2+\rho} (\log{P})^{2d+5}$, where $T_1 \ll \ldots T_d$ and $T_i \ll 1$ was used as before. Next we apply Corollary \ref{cor:averagewithresonance} with $k=2$ (again the assumptions guarantee that the quadratic form \eqref{def:Q_j} is indefinite) to find
 \begin{equation*}
  P^{\frac{1}{\beta_2}-\rho} q^{-1-\frac{1}{2 \beta_2}} \ll T_{2}^2 \ll q^2 P^{-2+\rho} (\log{P})^{2d+5}
 \end{equation*}
  and after rearranging
  \begin{equation}
    \label{eq:proof:couple_d-1}
    1 \ll P^{2\rho} (\log P)^{2d+5}u(P)^2 P^{\mathfrak p_2(d)},
  \end{equation}
  where
  \begin{equation*}
    \mathfrak{p}_2(d) := \mfrac{2}{(1+2\beta)} \big(3+\mfrac{1}{2 \beta_2}\big)-\big(2+\mfrac{1}{\beta_2}\big).
  \end{equation*}
  Considering again all cases in Table \ref{table:third} of the \hyperref[appendix-a]{Appendix A} shows that $\mathfrak p_2(d) <0$. Thus, inequality \eqref{eq:proof:couple_d-1} cannot hold if we increase $C_d >1$ and choose $\rho>0$ small enough. Again we conclude that the second alternative in Lemma \ref{lemma:dioph_approx} holds. The remaining steps are now the same as in the previous proof: We can apply Lemma \ref{lemma:BirDav:Lemma16} with $I= \{2,\ldots,d-1\}$.
\end{proof}

 By Corollary \ref{cor:first_coupling} we know that $S_2,\ldots ,S_d$ are coupled if $5 \leq d \leq 10$. Hence we may assume that $d \geq 11$ and then $S_2,\ldots, S_d$ are coupled as well if the signature of $Q$ is of the form $(d-1,1)$ or $(d-2,2)$. The last statement follows from (ii) of Corollary \ref{cor:first_coupling}, since $r \geq 2s$ holds for these cases. Thus, we have proven that $S_2,\ldots, S_d$ are coupled on $\mathcal{G}$, regardless of the signature $(r,s)$.

\begin{lemma}
 \label{lastcoupling}
  Under Assumption \ref{assumption:Theorem1} all exponential sums $S_{1},\ldots,S_{d}$ are coupled on $\calG$.
\end{lemma}

\begin{proof}
 By the previous discussion, we know that $S_2,\dots, S_d$ are coupled on $\mathcal{G}$. We can also assume that $d \ge 7$ and that the signature of $Q$ is not of the form $(d-1,1)$, since otherwise all coordinates are coupled, see Corollary \ref{cor:first_coupling}. Similar to the previous cases, we find
 \begin{equation}
   \label{eq:lastcoupling:proof}
   P^{\frac{3}{\beta_1}-3\rho} q^{-3-\frac{3}{2 \beta_1}} \ll T_{1}^3 \ll q^2 P^{-2+\rho} (\log{P})^{2d+5},
 \end{equation}
 where we removed the factor $T_d^4 T_1^1$ (by using $T_1 \ll \ldots T_d$ and $T_i \ll 1$) and applied Corollary \ref{cor:averagewithresonance} with $k=1$ (the assumptions are met since $Q$ is not of the form $(d-1,1)$). As before, the inequality \eqref{eq:lastcoupling:proof} can be rewritten as
 \begin{equation*}
    1 \ll P^{4 \rho} (\log P)^{2d+5}u(P)^2 P^{\mathfrak p_1(d)},
 \end{equation*}
 where
 \begin{equation*}
   \mathfrak{p}_1(d) := \mfrac{2}{(1+2\beta)} \big(5+\mfrac{3}{2 \beta_1}\big)-\big(2+\mfrac{3}{\beta_1}\big).
 \end{equation*}
 For every case, other than $\operatorname{sgn}(Q)=(\frac{d+3}{2},\frac{d-3}{2})$, we read off from Table \ref{table:third} in \hyperref[appendix-a]{Appendix A} that $\mathfrak p_1(d)<0$, thus yielding a contradiction. For $\text{sgn}(Q)=(\frac{d+3}{2},\frac{d-3}{2})$ and $2 \beta_1 = \frac{d+1}{d-5}$ we obtain also $\mathfrak p_1(d)= -\frac{6(d-5)}{d(d+1)}<0$. However, if $2 \beta_1 = \frac{d+3}{d-5}$, then $\mathfrak p_1(d) =0$. In this case the $(d-1)$-dimensional restriction of the quadratic form is of signature $(\frac{d+1}{2}+1,\frac{d-1}{2}-2)$ and hence we may remove one of the coordinates corresponding to $T_2,\dots,T_d$ to obtain a $(d-2)$-dimensional restriction of our quadratic form of signature $(\frac{d+1}{2},\frac{d-1}{2}-2)$. As in Corollary \ref{cor:averagewithresonance} (by applying Lemma \ref{averagewithresonance} to the aforementioned restriction of $Q$) we may deduce the inequality
 \begin{equation*}
    T_1T_l \gg P^{\frac{1}{\beta_2}-\rho} q^{-1-\frac{1}{2\beta_2}},
 \end{equation*}
 for some $2 \leq l \leq d$. Arguing again as above, we obtain
 \begin{equation*}
   P^{\frac{3}{\beta_2}-3\rho} q^{-3-\frac{3}{2\beta_2}} \ll T_{1}^3 \ll q^2 P^{-2+\rho} (\log{P})^{2d+5} u(P)^2 ,
 \end{equation*}
 which implies $1  \ll P^{-2+4\rho} (\log{P})^{2d+5} u(P)^2 P^{{\mathfrak p}_1(d)}$, where
 \begin{equation*}
    \mathfrak{p}_1(d) := \mfrac{2}{1+2 \beta} \big(5+\mfrac{3}{2\beta_2}\big)-\big(2+\mfrac{3}{\beta_2}\big)= -\mfrac{6(d-5)}{d(d+1)}<0.
 \end{equation*}
 We reach again a contradiction. Thus, the second alternative in Lemma \ref{lemma:dioph_approx} is valid. Since the previous considerations exhaust all cases, we can apply Lemma \ref{lemma:BirDav:Lemma16} with $I= \{1,\ldots,d-1\}$ and conclude that all coordinates are coupled on $\mathcal{G}$.
\end{proof}

\section{Proof of Theorem 1.3: Counting Approximants}
\label{sec:CoutingApproximants}

\noindent Finally, we are going to deduce a contradiction in form of an inconsistent inequality consisting of the lower bound for $N_j$, established in Corollary \ref{lemma:BirDav:Lemma12}, and the upper bound from Corollary \ref{cor:bound_fundamental} for the number of distinct pairs $(x,y)$.

\begin{proof}[Proof of Theorem \ref{MainTheorem}:]
 As shown in Subsection \ref{subsec:couple_all}, all coordinates can be coupled (under the Assumption \ref{assumption:Theorem1}) and therefore we can apply Corollary \ref{cor:bound_fundamental} with $k=0$ - in particular, we have $Q_k =Q$ - to find an upper bound for the number $N_j$ of all $(x_j,y_j)$: Since $x_1',y_1',\ldots,x_d',y_d'$ are determined as divisors of an $\alpha$-independent number $L \ll H^{10d}$, see Lemma \ref{lemma:BirDav:Lemma16}, Wigert's divisor bound (compare \eqref{eq:wigerts_bound}) implies that
 \begin{equation*}
  N_j^{2\beta} \ll H^{\frac{20d(d-1)}{\log \log H }} (\# \mathfrak{C}_0)^{2\beta} \ll H^{\frac{20d(d-1)}{\log \log H }} P^{-2} q^{2\beta+1} u(P)^{2\beta} (U_1 \ldots U_d)^\beta \big( \max_{i=1,\ldots,d} T_i^{-1} U_i^{-1/2}\big),
 \end{equation*}
 where we also used that $(4\beta+2)/d \le \beta$, which can be checked by considering the lower bound \eqref{eq:lowerbound:beta}. Next let $j \neq l$, where $l$ is an index for which the maximum of $T_i^{-1} U_i^{-1/2}$ is attained. Combined with the lower bound on $N_j$, obtained in Corollary \ref{lemma:BirDav:Lemma12}, we find
  \begin{equation}
    \label{eq:lastProof:eq1}
    (\log{P})^{- 4d \beta-1 } ( \mprod_{i=1}^d T_i )^{-2\beta} (T_j U_j^{\frac{1}{2}})^{2\beta} \ll H^\frac{20d(d-1)}{\log \log{H}} P^{-2} q^{2\beta+1} u(P)^{2\beta} (\mprod_{i=1}^d U_i)^\beta \, (T_l U_l^{\frac{1}{2}})^{-1}
  \end{equation}
  and this inequality can be simplified by using the notation
  \begin{equation*}
    V_i := U_i^{-1/2} T_i^{-1} (\log{P}).
  \end{equation*}
  Indeed, since $V_i \gg 1$ by \eqref{eq:UpperBoundForU_j}, we can rewrite \eqref{eq:lastProof:eq1} as
  \begin{align*}
    1 \ll (V_1 \ldots V_d)^{2\beta} \, V_j^{-2\beta} \, V_l^{-1} &\ll H^{-\frac{20d}{\log \log{H}}} u(P)^{2\beta} (\log{P})^{6d \beta+1}\\
    &\ll H^{-\frac{1}{\log \log H}} \leq \exp \big(-\tfrac{\log C_d}{\log \log C_d} \big),
  \end{align*}
  where $2\beta \geq 1$ was used. If $C_d \gg 1$ is chosen sufficiently large, we get a contradiction. Thus, our initial Assumption \ref{assumption:Theorem1} is false.
\end{proof}

\section*{Appendix A: Possible Signatures and Exponents}
\label{appendix-a}
\renewcommand{\theequation}{A.\arabic{equation}}
\noindent This Appendix constitutes sufficient preparation for the coupling argument: We determine all possible values of $\beta_k$ depending on the signature $(r,s)$ of $Q$ and give upper bounds for the exponents occuring in the iteration of the coupling argument.

\begin{center}
   \resizebox{15.70cm}{!}{
  \begin{tabular}{ | c | c || l | c || l | c || l | c | }
    \hline
    \multicolumn{8}{ |c |}{Even $d$}                                                                                                                                                                                                                                                                                                                                                                                                                                                                        \\
    \hline \hline
    $\text{Sign}(Q)$ \Tstrut \Bstrut                             & $2\beta$                                                   & $\text{Sign}(Q_{3})$                                     & $2\beta_{3}$                                                 & $\text{Sign}(Q_{2})$                                         & $2\beta_{2}$                                                   & $\text{Sign}(Q_{1})$                                       & $2\beta_{1}$                             \\  \hline \hline
    \multirow{5.4}{*}{$\big(\frac{d}{2},\frac{d}{2}\big)$}       & \multirow{5.6}{*}{$\frac{d+2}{d-4}$ \Tstrut\Bstrut}
                                                                 & $\big(\frac{d-6}{2},\frac{d}{2}\big)$ \Tstrut\Bstrut       & \multirow{5.6}{*}{$\frac{d}{d-6}$ \Tstrut\Bstrut}            & \multirow{2}{*}{$\big(\frac{d-4}{2},\frac{d}{2}\big)$}       &                                                                  & \multirow{3.7}{*}{$\big(\frac{d-2}{2},\frac{d}{2}\big)$}       &                                                                                                           \\
                                                                 &                                                            & $\big(\frac{d-4}{2},\frac{d-2}{2}\big)$ \Tstrut\Bstrut       &                                                              & \multirow{2.3}{*}{$\big(\frac{d-2}{2},\frac{d-2}{2}\big)$}       & \multirow{2.3}{*}{$\;\frac{d}{d-6}$ \Tstrut\Bstrut}            & \multirow{4.5}{*}{$\big(\frac{d}{2},\frac{d-2}{2}\big)$}       & \multirow{2.3}{*}{$\frac{d+2}{d-4}$}       \\
                                                                 &                                                            & $\big(\frac{d-2}{2},\frac{d-4}{2}\big)$ \Tstrut\Bstrut       &                                                              & \multirow{2.7}{*}{$\big(\frac{d}{2},\frac{d-4}{2}\big)$}         &                                                                &                                                                &                                          \\
                                                                 &                                                            & $\big(\frac{d}{2},\frac{d-6}{2}\big)$ \Tstrut\Bstrut         &                                                              &                                                                  &                                                                &                                                                &                                          \\  \hline \hline
    \multirow{5.4}{*}{$\big(\frac{d+2}{2},\frac{d-2}{2}\big)$}   & \multirow{5.6}{*}{$\frac{d+2}{d-4}$ \Tstrut\Bstrut}
                                                                 & $\big(\frac{d-4}{2},\frac{d-2}{2}\big)$ \Tstrut\Bstrut     & $\frac{d}{d-6}$                                              & \multirow{2}{*}{$\big(\frac{d-2}{2},\frac{d-2}{2}\big)$}     & \multirow{2}{*}{$\frac{d}{d-6}$}                                 & \multirow{3.7}{*}{$\big(\frac{d}{2},\frac{d-2}{2}\big)$}       & \multirow{5.4}{*}{$\frac{d+2}{d-4}$}                                                                        \\
                                                                 &                                                            & $\big(\frac{d-2}{2},\frac{d-4}{2}\big)$ \Tstrut\Bstrut       & $\frac{d}{d-6}$                                              & \multirow{2.3}{*}{$\big(\frac{d}{2},\frac{d-4}{2}\big)$}         & \multirow{2.3}{*}{$\; \frac{d}{d-6}$ \Tstrut\Bstrut}           & \multirow{4.5}{*}{$\big(\frac{d+2}{2},\frac{d-4}{2}\big)$}     &      \\
                                                                 &                                                            & $\big(\frac{d}{2},\frac{d-6}{2}\big)$ \Tstrut\Bstrut         & $\frac{d}{d-6}$                                              & \multirow{2.7}{*}{$\big(\frac{d+2}{2},\frac{d-6}{2}\big)$}       & \multirow{2.7}{*}{$\frac{d+2}{d-6}$}                           &                                                                &                                          \\
                                                                 &                                                            & $\big(\frac{d+2}{2},\frac{d-8}{2}\big)$ \Tstrut\Bstrut       & $\frac{d+2}{d-8}$                                            &                                                                  &                                                                &                                                                &                                          \\
    \hline \hline
    \multirow{5.4}{*}{$\big(\frac{d+4}{2},\frac{d-4}{2}\big)$}   & \multirow{5.6}{*}{$\frac{d+4}{d-4}$ \Tstrut\Bstrut}
                                                                 & $\big(\frac{d-2}{2},\frac{d-4}{2}\big)$ \Tstrut\Bstrut     & $\frac{d}{d-6}$                                              & \multirow{2}{*}{$\big(\frac{d}{2},\frac{d-4}{2}\big)$}       & \multirow{2}{*}{$\frac{d}{d-6}$}                                 & \multirow{3.7}{*}{$\big(\frac{d+2}{2},\frac{d-4}{2}\big)$}     & \multirow{3.7}{*}{$\frac{d+2}{d-4}$}                                                                      \\
                                                                 &                                                            & $\big(\frac{d}{2},\frac{d-6}{2}\big)$ \Tstrut\Bstrut         & $\frac{d}{d-6}$                                              & \multirow{2.3}{*}{$\big(\frac{d+2}{2},\frac{d-6}{2}\big)$}       & \multirow{2.3}{*}{$\; \frac{d+2}{d-6}$ \Tstrut\Bstrut}         & \multirow{4.5}{*}{$\big(\frac{d+4}{2},\frac{d-6}{2}\big)$}     & \multirow{4.5}{*}{$\frac{d+4}{d-6}$}     \\
                                                                 &                                                            & $\big(\frac{d+2}{2},\frac{d-8}{2}\big)$ \Tstrut\Bstrut       & $\frac{d+2}{d-8}$                                            & \multirow{2.7}{*}{$\big(\frac{d+4}{2},\frac{d-8}{2}\big)$}       & \multirow{2.7}{*}{$\frac{d+4}{d-8}$}                           &                                                                &                                          \\
                                                                 &                                                            & $\big(\frac{d+4}{2},\frac{d-10}{2}\big)$ \Tstrut\Bstrut      & $\frac{d+4}{d-10}$                                           &                                                                  &                                                                &                                                                &                                          \\
    \hline \hline
    \multirow{4.5}{*}{$\big(\frac{d+2l}{2},\frac{d-2l}{2}\big)$} & \multirow{5.6}{*}{$\frac{d+2l}{d-2l}$ \Tstrut\Bstrut}
                                                                 & $\big(\frac{d+2l-6}{2},\frac{d-2l}{2}\big)$ \Tstrut\Bstrut & $\frac{d+2l-6}{d-2l}$                                        & \multirow{2}{*}{$\big(\frac{d+2l-4}{2},\frac{d-2l}{2}\big)$} & \multirow{2}{*}{$\frac{d+2l-4}{d-2l}$}                           & \multirow{3.7}{*}{$\big(\frac{d+2l-2}{2},\frac{d-2l}{2}\big)$} & \multirow{3.7}{*}{$\frac{d+2l-2}{d-2l}$}                                                                  \\
    \multirow{4.5}{*}{$l \geq 3$}                                &                                                            & $\big(\frac{d+2l-4}{2},\frac{d-2l-2}{2}\big)$ \Tstrut\Bstrut & $\frac{d+2l-4}{d-2l-2}$                                      & \multirow{2.3}{*}{$\big(\frac{d+2l-2}{2},\frac{d-2l-2}{2}\big)$} & \multirow{2.3}{*}{$\; \frac{d+2l-2}{d-2l-2}$ \Tstrut\Bstrut}   & \multirow{4.5}{*}{$\big(\frac{d+2l}{2},\frac{d-2l-2}{2}\big)$} & \multirow{4.5}{*}{$\frac{d+2l}{d-2l-2}$} \\
                                                                 &                                                            & $\big(\frac{d+2l-2}{2},\frac{d-2l-4}{2}\big)$ \Tstrut\Bstrut & $\frac{d+2l-2}{d-2l-4}$                                      & \multirow{2.7}{*}{$\big(\frac{d+2l}{2},\frac{d-2l-4}{2}\big)$}   & \multirow{2.7}{*}{$\frac{d+2l}{d-2l-4}$}                       &                                                                &                                          \\
                                                                 &                                                            & $\big(\frac{d+2l}{2},\frac{d-2l-6}{2}\big)$ \Tstrut\Bstrut   & $\frac{d+2l}{d-2l-6}$                                        &                                                                  &                                                                &                                                                &                                          \\
    \hline
  \end{tabular}
  } 
\end{center}
\begin{center}
   \resizebox{15.70cm}{!}{
  \begin{tabular}{ | c | c || l | c || l | c || l | c | }
    \hline
    \multicolumn{8}{ |c |}{Odd $d$}                                                                                                                                                                                                                                                                                                                                                                                                                                                                                       \\
    \hline \hline
    $\text{Sign}(Q)$ \Tstrut \Bstrut                                 & $2\beta$                                                     & $\text{Sign}(Q_{3})$                                     & $2\beta_{3}$                                                   & $\text{Sign}(Q_{2})$                                         & $2\beta_{2}$                                                     & $\text{Sign}(Q_{1})$                                         & $2\beta_{1}$                               \\  \hline \hline
    \multirow{5.4}{*}{$\big(\frac{d+1}{2},\frac{d-1}{2}\big)$}       & \multirow{5.6}{*}{$\frac{d+3}{d-3}$ \Tstrut\Bstrut}
                                                                     & $\big(\frac{d-5}{2},\frac{d-1}{2}\big)$ \Tstrut\Bstrut         & $\frac{d-1}{d-7}$                                            & \multirow{2}{*}{$\big(\frac{d-3}{2},\frac{d-1}{2}\big)$}       &                                                                  & \multirow{3.7}{*}{$\big(\frac{d-1}{2},\frac{d-1}{2}\big)$}       &                                                                                                               \\
                                                                     &                                                              & $\big(\frac{d-3}{2},\frac{d-3}{2}\big)$ \Tstrut\Bstrut       & $\frac{d-1}{d-7}$                                              & \multirow{2.3}{*}{$\big(\frac{d-1}{2},\frac{d-3}{2}\big)$}       & \multirow{2.3}{*}{$\;\frac{d+1}{d-5}$ \Tstrut\Bstrut}            & \multirow{4.5}{*}{$\big(\frac{d+1}{2},\frac{d-3}{2}\big)$}       & \multirow{2.3}{*}{$\frac{d+1}{d-5}$}       \\
                                                                     &                                                              & $\big(\frac{d-1}{2},\frac{d-5}{2}\big)$ \Tstrut\Bstrut       & $\frac{d-1}{d-7}$                                              & \multirow{2.7}{*}{$\big(\frac{d+1}{2},\frac{d-5}{2}\big)$}       &                                                                  &                                                                  &                                            \\
                                                                     &                                                              & $\big(\frac{d+1}{2},\frac{d-7}{2}\big)$ \Tstrut\Bstrut       & $\frac{d+1}{d-7}$                                              &                                                                  &                                                                  &                                                                  &                                            \\  \hline \hline
    \multirow{5.4}{*}{$\big(\frac{d+3}{2},\frac{d-3}{2}\big)$}       & \multirow{5.6}{*}{$\frac{d+3}{d-3}$ \Tstrut\Bstrut}
                                                                     & $\big(\frac{d-3}{2},\frac{d-3}{2}\big)$ \Tstrut\Bstrut       & $\frac{d-1}{d-7}$                                            & \multirow{2}{*}{$\big(\frac{d-1}{2},\frac{d-3}{2}\big)$}       & \multirow{2}{*}{$\frac{d+1}{d-5}$}                               & \multirow{3.7}{*}{$\big(\frac{d+1}{2},\frac{d-3}{2}\big)$}       & \multirow{3.7}{*}{$\frac{d+1}{d-5}$}                                                                          \\
                                                                     &                                                              & $\big(\frac{d-1}{2},\frac{d-5}{2}\big)$ \Tstrut\Bstrut       & $\frac{d-1}{d-7}$                                              & \multirow{2.3}{*}{$\big(\frac{d+1}{2},\frac{d-5}{2}\big)$}       & \multirow{2.3}{*}{$\; \frac{d+1}{d-5}$ \Tstrut\Bstrut}           & \multirow{4.5}{*}{$\big(\frac{d+3}{2},\frac{d-5}{2}\big)$}       & \multirow{4.5}{*}{$\frac{d+3}{d-5}$}       \\
                                                                     &                                                              & $\big(\frac{d+1}{2},\frac{d-7}{2}\big)$ \Tstrut\Bstrut       & $\frac{d+1}{d-7}$                                              & \multirow{2.7}{*}{$\big(\frac{d+3}{2},\frac{d-7}{2}\big)$}       & \multirow{2.7}{*}{$\frac{d+3}{d-7}$}                             &                                                                  &                                            \\
                                                                     &                                                              & $\big(\frac{d+3}{2},\frac{d-9}{2}\big)$ \Tstrut\Bstrut       & $\frac{d+3}{d-9}$                                              &                                                                  &                                                                  &                                                                  &                                            \\
    \hline \hline
    \multirow{5.4}{*}{$\big(\frac{d+5}{2},\frac{d-5}{2}\big)$}       & \multirow{5.6}{*}{$\frac{d+5}{d-5}$ \Tstrut\Bstrut}
                                                                     & $\big(\frac{d-1}{2},\frac{d-5}{2}\big)$ \Tstrut\Bstrut       & $\frac{d-1}{d-7}$                                            & \multirow{2}{*}{$\big(\frac{d+1}{2},\frac{d-5}{2}\big)$}       & \multirow{2}{*}{$\frac{d+1}{d-5}$}                               & \multirow{3.7}{*}{$\big(\frac{d+3}{2},\frac{d-5}{2}\big)$}       & \multirow{3.7}{*}{$\frac{d+3}{d-5}$}                                                                          \\
                                                                     &                                                              & $\big(\frac{d+1}{2},\frac{d-7}{2}\big)$ \Tstrut\Bstrut       & $\frac{d+1}{d-7}$                                              & \multirow{2.3}{*}{$\big(\frac{d+3}{2},\frac{d-7}{2}\big)$}       & \multirow{2.3}{*}{$\; \frac{d+3}{d-7}$ \Tstrut\Bstrut}           & \multirow{4.5}{*}{$\big(\frac{d+5}{2},\frac{d-7}{2}\big)$}       & \multirow{4.5}{*}{$\frac{d+5}{d-7}$}       \\
                                                                     &                                                              & $\big(\frac{d+3}{2},\frac{d-9}{2}\big)$ \Tstrut\Bstrut       & $\frac{d+3}{d-9}$                                              & \multirow{2.7}{*}{$\big(\frac{d+5}{2},\frac{d-9}{2}\big)$}       & \multirow{2.7}{*}{$\frac{d+5}{d-9}$}                             &                                                                  &                                            \\
                                                                     &                                                              & $\big(\frac{d+5}{2},\frac{d-11}{2}\big)$ \Tstrut\Bstrut      & $\frac{d+5}{d-11}$                                             &                                                                  &                                                                  &                                                                  &                                            \\
    \hline \hline
    \multirow{4.5}{*}{$\big(\frac{d+2l+1}{2},\frac{d-2l-1}{2}\big)$} & \multirow{5.6}{*}{$\frac{d+2l+1}{d-2l-1}$ \Tstrut\Bstrut}
                                                                     & $\big(\frac{d+2l-5}{2},\frac{d-2l-1}{2}\big)$ \Tstrut\Bstrut & $\frac{d+2l-5}{d-2l-1}$                                      & \multirow{2}{*}{$\big(\frac{d+2l-3}{2},\frac{d-2l-1}{2}\big)$} & \multirow{2}{*}{$\frac{d+2l-3}{d-2l-1}$}                         & \multirow{3.7}{*}{$\big(\frac{d+2l-1}{2},\frac{d-2l-1}{2}\big)$} & \multirow{3.7}{*}{$\frac{d+2l-1}{d-2l-1}$}                                                                    \\
    \multirow{4.5}{*}{$l \geq 3$}                                    &                                                              & $\big(\frac{d+2l-3}{2},\frac{d-2l-3}{2}\big)$ \Tstrut\Bstrut & $\frac{d+2l-3}{d-2l-3}$                                        & \multirow{2.3}{*}{$\big(\frac{d+2l-1}{2},\frac{d-2l-3}{2}\big)$} & \multirow{2.3}{*}{$\; \frac{d+2l-1}{d-2l-3}$ \Tstrut\Bstrut}     & \multirow{4.5}{*}{$\big(\frac{d+2l+1}{2},\frac{d-2l-3}{2}\big)$} & \multirow{4.5}{*}{$\frac{d+2l+1}{d-2l-3}$} \\
                                                                     &                                                              & $\big(\frac{d+2l-1}{2},\frac{d-2l-5}{2}\big)$ \Tstrut\Bstrut & $\frac{d+2l-1}{d-2l-5}$                                        & \multirow{2.7}{*}{$\big(\frac{d+2l+1}{2},\frac{d-2l-5}{2}\big)$} & \multirow{2.7}{*}{$\frac{d+2l+1}{d-2l-5}$}                       &                                                                  &                                            \\
                                                                     &                                                              & $\big(\frac{d+2l+1}{2},\frac{d-2l-7}{2}\big)$ \Tstrut\Bstrut & $\frac{d+2l+1}{d-2l-7}$                                        &                                                                  &                                                                  &                                                                  &                                            \\
    \hline
  \end{tabular}
  }
\end{center}

Note that in both tables the last case in every row is the worst when compared to $\beta$. Thus, considering theses cases, one can derive the following bound on the exponent $\mathfrak p_i(d)$ appearing in the iteration of the coupling argument, see Lemmas \ref{firstcoupling} - \ref{lastcoupling}.

  \begin{center}
    \captionsetup{type=table}
    \begin{tabular}{ | c | c | c | c |}
      \hline
      $\text{Sign}(Q)$ \Tstrut \Bstrut                                          & \hspace{0.4cm} $\mathfrak p_3(d) \leq$ \hspace{0.4cm} & \hspace{0.4cm} $\mathfrak p_2(d) \leq$ \hspace{0.4cm} & \hspace{0.4cm} $\mathfrak p_1(d) \leq$ \hspace{0.4cm} \\  \hline \hline
      $\big(\frac{d}{2},\frac{d}{2}\big)$ \Tstrut\Bstrut                        & $-\frac{6d-4}{d(d-1)}$ \Tstrut\Bstrut
                                                                                & $-\frac{6(d-2)}{d(d-1)}$                                 &    $-\frac{6}{d-1}$                                            \\  \hline
      $\big(\frac{d+2}{2},\frac{d-2}{2}\big)$ \Tstrut\Bstrut                    & $-\frac{14}{3(d-1)}$ \Tstrut\Bstrut
                                                                                & $-\frac{4}{d-1}$                                         &                  $-\frac{6}{d-1}$                              \\  \hline
      $\big(\frac{d+2l}{2},\frac{d-2l}{2}\big)$,  $l \geq 2$ \Tstrut\Bstrut     & $-\frac{2(2l-1)}{d}$ \Tstrut\Bstrut
                                                                                & $-\frac{4(l-1)}{d}$                                   &                             $-\frac{2(2l-3)}{d}$                   \\  \hline \hline
      $\big(\frac{d+1}{2},\frac{d-1}{2}\big)$ \Tstrut\Bstrut                    & $-\frac{16}{3(d+1)}$ \Tstrut\Bstrut
                                                                                & $-\frac{6(d-1)}{d(d+1)}$                             &                         $-\frac{6(d-5)}{d(d+1)}$                       \\  \hline
      $\big(\frac{d+3}{2},\frac{d-3}{2}\big)$ \Tstrut\Bstrut                    & $-\frac{4}{d}$ \Tstrut\Bstrut
                                                                                & $-\frac{2}{d}$                                         &                          $*$                      \\  \hline
     \ \ $\big(\frac{d+2l+1}{2},\frac{d-2l-1}{2}\big)$,  $l \geq 2$ \ \ \Tstrut\Bstrut & $-\frac{4l}{d}$ \Tstrut\Bstrut
                                                                                & $-\frac{2(2l-1)}{d}$                                &                           $-\frac{4(l-1)}{d}$                     \\  \hline
    \end{tabular}
    \captionof{table}{Bounds on the exponents $\mathfrak p_1(d),\mathfrak p_2(d),\mathfrak p_3(d)$}\label{table:third}
  \end{center}

\section*{Appendix B: Kernels with fast-decaying Fourier transforms}
\renewcommand{\theequation}{B.\arabic{equation}}
\label{appendix-b}
\noindent In this Appendix we give a complete proof of Lemma \ref{lemma:Improvement:BirDav:Lemma1} showing the existence of compactly-supported kernels with fast-decaying Fourier transforms. The proof given here is elementary and based on arguments presented in \cite{Bhattacharya-RangaRao:2010} (see Theorem 10.2).

\begin{proof}[Proof of Lemma \ref{lemma:Improvement:BirDav:Lemma1}]
 During this proof we write
 \begin{equation*}
   U([-a,a]) = (2a)^{-1} 1_{[-a,a]}
 \end{equation*}
 for the density function of the uniform distribution on some interval $[-a,a]$, $a>0$, whose Fourier transform is given by
 \begin{equation}
   \label{app_b:fourier}
   \widehat{U}([-a,a])(t) = \frac{\sin(2\pi at)}{2\pi at}.
 \end{equation}
 Based on this simple kernel we will construct an infinite convolution product: First we shall make use of condition \eqref{eq:ingham:condition}, that is
 $\int_1^\infty \frac{1}{\alpha u(\alpha)} \, \dif \alpha < \infty$, by noting that there exists an integer $n_0 \in \mathbb{N}$ and a non-decreasing sequence of non-negative numbers $(a_n)_{n \in \mathbb{N}}$ given by
 \begin{equation*}
   a_n = \begin{cases}
             \frac{\mathrm{e}}{n_0 u(n_0) } & \text{ if } 1 \leq n  \le n_0 \\ \frac{\mathrm{e}}{n u(n)} & \text{ if } n > n_0
          \end{cases}
  \end{equation*}
  such that
   \begin{equation}
     \label{app_b:eq0}
     \sum_{n=1}^\infty a_r = \frac{\mathrm{e}}{u(n_0)} + \mathrm{e} \sum_{n=n_0+1}^\infty \frac{1}{n u(n)} \leq 1.
   \end{equation} 
   It remains to check that the sequence
   \begin{equation*}
     \psi_n := U([-a_1,a_1]) \ast \ldots \ast U([-a_n,a_n])
   \end{equation*}
   is uniformly convergent and satisfies the properties claimed in Lemma \ref{lemma:Improvement:BirDav:Lemma1}. To do this, we first verify that any $\psi_n$, $n \geq 2$, is Lipschitz continuous with Lipschitz constant $ 1/( 4 a_1 a_2)$. In fact, if $0 <b \le  a$, a simple calculation shows
   \begin{equation}
     \label{app_b:eq1}
     U([-a,a]) \ast U([-b,b])(t) = \begin{cases}
                                     0            & \text{if } \abs{t} \ge a+b \\
                                     \frac{1}{2a} & \text{if } \abs{t} \le a-b \\
                                     \frac{a+b -\abs{t}}{4ab} & \text{else}
                                   \end{cases}.
   \end{equation}
   Thus the above remark is true for $n=2$. The general case follows by induction:
   \begin{equation*}
     \abs{u_{n+1}(s)-u_{n+1}(t)} \leq \frac{1}{2 a_{n+1} } \int_{-a_{n+1}}^{a_{n+1}} \abs{u_n(s-h)-u_n(t-h)} \, \dif h \leq \frac{1}{4 a_1 a_2} \abs{t-s}.
   \end{equation*}
   Proceeding in the same manner, we see for any $n \geq 1$ that
   \begin{align*}
     \abs{u_{n+1}(t) - u_n(t)} &\leq \frac{1}{2a_{n+1}} \int_{-a_{n+1}}^{a_{n+1}} \abs{u_n(t-h)-u_n(t)} \, \dif h \\
     &\leq  \int_{-a_{n+1}}^{a_{n+1}} \frac{\abs{h}}{8 a_1 a_2 a_{n+1}} \, \dif h = \frac{a_{n+1}}{8 a_1 a_2}.
   \end{align*}
   In view of \eqref{app_b:eq0} this shows that $(\psi_n)_{n \in \mathbb{N}}$ is uniformly convergent, say to $\psi$, and thus a continuous probability density. $\psi$ is also symmetric, since by construction any $\psi_n$ is symmetric, and supported in $[-1,1]$, since $\psi_n$ has compact support lying in $[-\sum_{k=1}^n a_k, \sum_{k=1}^n a_k] \subset [-1,1]$, compare \eqref{app_b:eq0}. By induction we also find that $\psi_n$ is a $\mathcal{C}^{(n-2)}$-function (if $n \geq 1$) with  
   \begin{equation}
     \label{app_b:eq2}
     \psi_{n+1}^{(k+1)}(t) = \frac{1}{2 a_{n+1}} \Big( \psi_n^{(k)}(t+a_{n+1}) -  \psi_n^{(k)} (t-a_{n+1})\Big)
   \end{equation}
   in the range $0 \leq k \leq n-2$. The last line implies (again by induction) that $\psi_{k+2}^{(k)}$ is Lipschitz continuous with growing Lipschitz constant $L_k$ given inductively by $L_{k+1} = L_{k}/a_{k+3}$ and $L_0= 1/(4a_1 a_2)$. This in turn shows that $\abs{\psi^{(k)}_{n+1}(t) - \psi^{(k)}_n(t)} \leq L_k a_{n+1}$ for all $n \geq k +2$. Thus, we have confirmed the uniform convergence of any derivative, i.e.\ $\psi$ is smooth.
   
   Next we prove part (ii) of Lemma \ref{lemma:Improvement:BirDav:Lemma1} by induction: For $n=1$ the statement is trivial and for $n=2$ this follows at once from \eqref{app_b:eq1}. If $n \ge 3$ then we have 
   \begin{equation*}
     \psi_{n+1}'(t) = \frac{1}{2a_{n+1}} \{ \psi_n(t+a_{n+1}) - \psi_n(t-a_{n+1})\},
   \end{equation*}
   that is a special case of \eqref{app_b:eq2}. At this point we may use the symmetry of $\psi_n$ in order to conclude that both $\psi_{n+1}'(t) \ge 0$ if $t \le 0$ and $\psi_{n+1}'(t) \le 0$ if $t \ge 0$ hold, as claimed. Letting $n\rightarrow \infty$ yields part (ii) for $\psi$. In particular, $\psi$ has a global maximum in $t=0$ and it follows that $2 \psi(0) \geq \int \psi(t) \, \dif t =1$, i.e.\ the second part of (i) holds as well.
   
   Finally, it remains to prove part (iii) of Lemma \ref{lemma:Improvement:BirDav:Lemma1}. The uniform convergence combined with the explicit formula \eqref{app_b:fourier} implies the representation
   \begin{equation*}
     \widehat{\psi}(t) = \prod_{n=1}^\infty \Big( \frac{\sin(2\pi a_n t)}{2\pi a_n t} \Big)
   \end{equation*}
   as an infinite product with uniform convergence on compact sets. Note that \eqref{eq:ingham:condition} necessarily implies $u(t) \rightarrow \infty$ if $t \rightarrow \infty$ and therefore there exists a $t_0 >0$ such that $u(t) \ge 1$ for all $t \ge t_0$. For any $\abs{t} \ge t_0$ we have the bound
   \begin{equation*}
      \abs{\widehat{\psi}(t)} \leq \prod_{k=1}^n \Big( \frac{1}{2 \pi \abs{a_k t}} \Big) \leq  \frac{1}{\abs{a_n t}^n} = \Big( \frac{n u(n) }{ \mathrm{e} \abs{t}} \Big)^n.
   \end{equation*}
   Thus, taking $n = \lfloor \abs{t} u( \abs{t})^{-1} \rfloor$ (i.e.\ the integer part of $\abs{t}u(\abs{t})^{-1}$) yields
   \begin{equation*}
         \abs{\widehat{\psi}(t)} \leq \Big( \frac{u(n)}{\mathrm{e} \abs{t}} \Big)^n \le \mathrm{e}^{-n} \ll \exp \{ - \abs{t}u(\abs{t})^{-1} \}.
   \end{equation*}
   In the last line we used that $u$ is non-decreasing and that $\abs{t} \ge n$, since $\abs{t} \ge t_0$. This completes the proof of Lemma \ref{lemma:Improvement:BirDav:Lemma1}.
\end{proof}

\section*{References}
\nocite{*}
\renewcommand*{\bibfont}{\small}
\printbibliography[heading=none]

\end{document}